\input ./style/arxiv-general.cfg
\documentclass[MSNbibl,number,citesort,seceqn,dvips]{arxbj}
\makeatletter
   \@ifpackageloaded{graphicx}{}{\usepackage{graphicx}}
\makeatother
\usepackage{upgreek}

\aid{0}
\volume{21}
\issue{4}
\pubyear{2015}
\firstpage{2242}
\lastpage{2288}
\doi{10.3150/14-BEJ643} 
\docsubty{FLA}

\makeatletter

\newcommand{\rrvert}{\vert}
\newcommand{\rrVert}{\Vert}
\newcommand{\llvert}{\vert}
\newcommand{\llVert}{\Vert}
\newtheorem{theorem}{Theorem}[section]
\newtheorem{lemma}[theorem]{Lemma}
\newtheorem{proposition}[theorem]{Proposition}
\newtheorem{cor}[theorem]{Corollary}
\newremark{remark}[theorem]{Remark}
\newproclaim{assume}[theorem]{Assumption}

\newcommand{\p}{\mathbb{P}}
\newcommand{\sign}{\mathrm{sign}}

\newcommand{\R}{\mathbb{R}}
\renewcommand{\tilde}{\widetilde}

\newcommand{\cF}{{{\mathcal F}}}
\newcommand{\cG}{{{\mathcal G}}}
\newcommand{\cK}{{{\mathcal K}}}
\newcommand{\cL}{{{\mathcal L}}}
\newcommand{\cT}{{{\mathcal T}}}

\newcommand{\dd}{\,\mathrm{d}}
\newcommand{\ddd}{\mathrm{d}}

\newcommand{\limtwo}[2]{\mathop{\lim_{#1}}_{#2}} 

\newcommand{\bbC}{{{\mathbb C}} }
\newcommand{\bbE}{{{\mathbb E}} }
\newcommand{\bbL}{{{\mathbb L}} }
\newcommand{\bbP}{{{\mathbb P}} }
\newcommand{\bbR}{{{\mathbb R}} }

\newcommand{\ga}{\alpha}
\newcommand{\gb}{\beta}
\newcommand{\gd}{\delta}
\newcommand{\gep}{\varepsilon} 
\newcommand{\gp}{\varphi}
\newcommand{\gr}{\rho}
\newcommand{\gO}{\Omega}
\newcommand{\gl}{\lambda}
\newcommand{\gs}{\sigma}

\newcommand{\spect}{\mathtt{Sp}}

\newcommand{\dens}{\mathtt{f}}

\newcommand{\e}{_{\varepsilon}}
\newcommand{\me}{_{\mu,\varepsilon}}
\newcommand{\ee}{\mathbb{E}}

\newcommand{\eqref}[1]{(\ref{#1})}

\def\sfrac#1#2{#1/#2}

\def\afrac#1#2{#1/(#2)}
\def\vafrac#1#2{(#1)/(#2)}

\makeatother

\begin{document}
\begin{frontmatter}

\title{Weak noise and non-hyperbolic unstable fixed points: Sharp estimates on transit and~exit~times}
\runtitle{Weak noise and non-hyperbolic unstable fixed points}

\begin{aug}
\author[A]{\inits{G.}\fnms{Giambattista}~\snm{Giacomin}} \and
\author[A]{\inits{M.}\fnms{Mathieu}~\snm{Merle}\corref{}\ead[label=e1]{merle@math.univ-paris-diderot.fr}}
\address[A]{Universit\'{e} Paris Diderot and Laboratoire de Probabilit\'{e}s et Mod\`{e}les Al\'{e}atoires (CNRS),
U.F.R. Math\'{e}matiques, Case 7012, B\^{a}t. Sophie Germain, 75205 Paris C\'{e}dex 13, France }
\end{aug}

\received{\smonth{7} \syear{2013}}
\revised{\smonth{5} \syear{2014}}

%
\begin{abstract}
We consider certain one dimensional ordinary stochastic differential
equations driven by additive Brownian motion
of variance $\gep^2$. When $\gep=0$ such equations have an
unstable non-hyperbolic fixed point
and the drift near such a point
has a power law behavior. For $\gep>0$ small, the fixed point property
disappears, but it is replaced by a random escape or
transit time which diverges as $\gep\searrow0$. We show that this
random time, under
suitable (easily guessed) rescaling, converges to a limit random
variable that essentially depends only
on the power exponent associated to the fixed point. Such random
variables, or laws, have therefore a universal character
and they arise of course in a variety of contexts.
We then obtain quantitative sharp estimates, notably tail properties,
on these
universal laws.
\end{abstract}

%
\end{frontmatter}

\section{Introduction}
\label{sec:intro}
\subsection{Effect of noise on non-hyperbolic unstable points}
Noise perturbations on dynamical systems lead to a variety of phenomena
and many are of crucial interest in understanding the dynamics of real
systems \cite{cf:HL,cf:exc-elem}. Here we focus on a basic issue that
has been repeatedly addressed in various domains \cite{cf:arecchi,cf:CFB,cf:LLB,cf:SH}: the effect of noise
on stationary non-hyperbolic points of one dimensional Ordinary
Differential Equations (ODE). The basic question we have in mind
is easily stated at an informal level:
consider the Stochastic Differential Equation (SDE)
%
\begin{equation}
\label{eq:SDE0} \mathrm{d} X_t = -U'(X_t)
\dd t + \gep\dd W_t ,
\end{equation}
where $\gep\ge0$, $\{W_t\}_{t\ge0}$ is a standard Brownian motion,
and $U(\cdot)$ is a smooth function such that $U'(0)=U''(0)=0$, that
is $0$ is a non-hyperbolic fixed point for the case
$\gep=0$. We also require such fixed point to be unstable: the cases
we have in mind are for example
%
\begin{equation}
\label{eq:example3} U(x) = -\frac{x^3}6 \quad \mbox{and}\quad  U(x) = \sin(x)-x ,
\end{equation}
that is the case in which $0$ is a saddle-node fixed point,
and
%
\begin{equation}
\label{eq:example4} U(x) = -\frac{x^4}4 \quad \mbox{and}\quad  U(x) = - \bigl( 1-\cos(x)
\bigr)^2 ,
\end{equation}
that is the case is which $0$ is a symmetric non-hyperbolic unstable
point: we may focus on these example for the sake of informal
discussion and we refer to them as the cases $d=3$ and $d=4$, in
conformity with the rest of the paper in
which we will address the case in which $U(x)$ is \textsl{roughly}
proportional to $x^d$ in a neighborhood of zero.

Switching from $\gep=0$ to $\gep>0$ will have the obvious drastic
effect on the solution. Nonetheless, if $\gep$ is small, it will
require a long time to leave the neighborhood of the origin. Two
comments are in order:
\begin{enumerate}[(2)]
\item[(1)] in general \eqref{eq:SDE0} does not admit a global (strong)
solution: in fact for the first case in both \eqref{eq:example3}
and \eqref{eq:example4} for $\gep>0$ and any choice of $X_0$ the
solution to \eqref{eq:SDE0} has a finite explosion time (see, e.g., \cite{cf:KaSh}, Section 5.5.C, or \cite{cf:groisman});
\item[(2)] the cases of $d=3$ and $4$ or, more generally, $d$ odd or even,
are different and in the former case we will
be interested in $X_0\in[-\infty, 0)$ so that $\lim_{t \to
\infty} X_t =0$ for $\gep=0$, while in the even $d$ case the most interesting
choice is $X_0=0$.
\end{enumerate}

We are after understanding the distribution of the time for going
through the saddle point (for $d$ odd) and
the distribution of the time of escape from $0$ (for $d$ even), in the
small $\gep$ limit. More precisely:
\begin{itemize}
\item In the cases in \eqref{eq:example3}, consider the first hitting
time $\tau_{a, \gep}$ of $a \in(0, \infty]$ or $a \in(0, \uppi)$
(according to whether we consider the first or second case) for $X_0\in
[-\infty,0)$ or $X_0\in(-\uppi,0)$.
It is easy to see that $\lim_{\gep\to0} \tau_{a, \gep}=+ \infty$
and by scaling argument is not difficult to guess that $ \tau_{a, \gep
}\approx\gep^{-2(d-2)/d}$ (e.g., \cite{cf:CFB,cf:LLB,cf:SH}, the
argument is also given explicitly in Section~\ref{sec:cubic})
and it is natural to expect that $\gep^{2(d-2)/d} \tau_{a, \gep}$
converges in law as $\gep\searrow0$ to a random variable $T_d$ that
does not depend on $a$, nor on $X_0$, nor on whether we have chosen the
first or second example. This has been argued for example in \cite{cf:CFB,cf:LLB,cf:SH} where
this claim is substantiated by explicit computations of mean and
variance of $T_d$ and by numerical computations.
%
\item In the cases in \eqref{eq:example4}, consider for $X_0=0$ the
first hitting time $\tau_{a, \gep}$, of $\pm a$ with $a$ either in
$(0, \infty]$ or in $(0, \uppi)$ according to which of the cases we consider.
Once again scaling arguments \cite{cf:arecchi,cf:pitchfork} suggest
that $\gep^{2(d-2)/d} \tau_{a, \gep}$ converges in law as $\gep
\searrow0$ to a random variable $T_d$ that
does not depend on $a$ nor on which of the two cases we have chosen.
\end{itemize}

\begin{figure}

\includegraphics{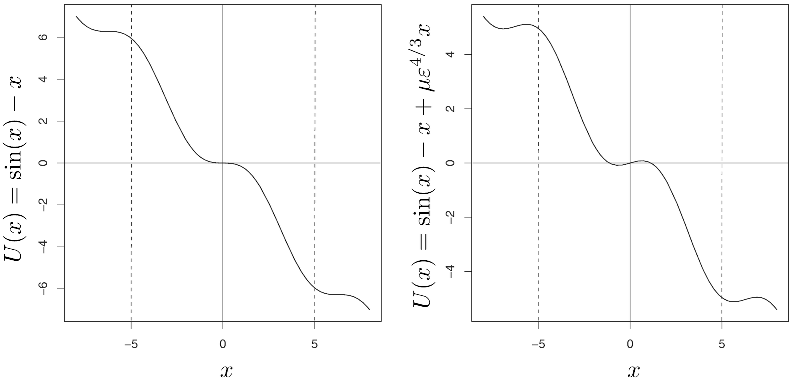}

\caption{A case with $d=3$. The potential on the left has saddle
points at integer multiples of $2 \uppi$. We just focus on the saddle
point at $0$ and we aim at capturing the leading asymptotic behavior of
the hitting time of $5$, or any other point
in $(0, 2 \uppi)$, starting from $-5$, or any other point
in $(-2 \uppi,0)$. For $\gep$ small, but positive, the effect of the
noise if essentially negligible except in a small neighborhood of $0$, which
turns out to be $\mathrm{O}(\gep^{2/3})$: in this neighborhood the noise
eventually drives $X_t$ through the saddle point, but this transit
event takes a long time (of order $\gep^{-2/3}$).
The phenomenon does not change qualitatively for suitable $\gep$-small
perturbations (figure on the right), but note that
the scaling in $\gep$ of the term containing $\mu$ has been carefully
chosen so that also this term contributes
to determine the limit transit time that in fact will be a new random
variable $T_{3,\mu}$, see Section~\protect\ref{sec:main}.}\label{fig:1}
\end{figure}

\begin{figure}

\includegraphics{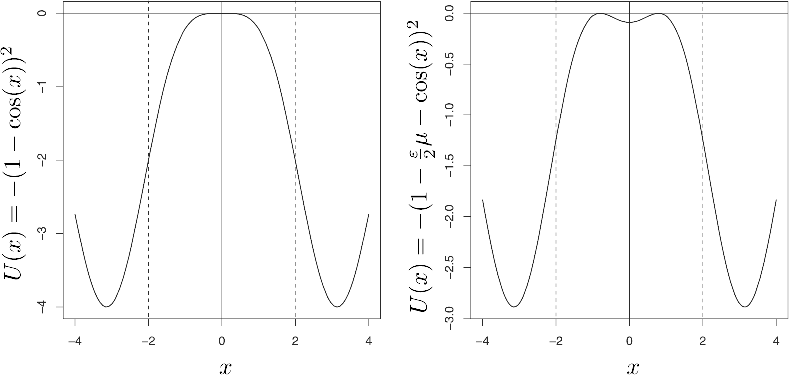}

\caption{A case with $d=4$. The potential on the left has unstable
non-hyperbolic points at integer multiples of $2 \uppi$. We just focus on
the unstable point at $0$ and we aim at capturing the leading
asymptotic behavior of the hitting time of $\{-2,2\}$, starting from
$0$. For $\gep$ small, but positive, the effect of the noise if
essentially negligible except in a small neighborhood of $0$, which
turns out to be $\mathrm{O}(\gep^{1/2})$: in this neighborhood the noise
eventually drives $X_t$ away from $0$: for this specific
potential $X_t$ will reach a neighborhood of $\uppi$ or $-\uppi$ and will
stay there for a time which to leading order is an exponential random
variable with mean $\exp(c \gep^{-1})$, for some $c>0$ that can be
computed by Large Deviations technics \cite{cf:FW}.
The escape from $0$ happens instead on a much shorter time scale ($\gep^{-1}$).
The phenomenon does not change qualitatively for suitable $\gep$-small
perturbations (figure on the right),
but the very same considerations at the end of the caption of Figure~\protect\ref{fig:1} apply in this case too.}
\label{fig:2}
\end{figure}

The applied literature based on \eqref{eq:SDE0} with the type of
potentials we are looking at is extremely vast, and the focus
on understanding $T_d$ is often at the heart of the analyses. We
mention here for example
the relevance of the odd $d$ case (saddle-node) in the context of
modeling excitable systems \cite{cf:exc-elem,cf:LLB,cf:SH} and
in this context, it is very natural to consider the extension to the
\textsl{weakly tilted case} of the left inset in Figure~\ref{fig:1}:
the saddle node can be in fact viewed as the critical or marginal case
of a saddle-node bifurcation
(this is
going to be taken up in Section~\ref{sec:gen}). We signal also the recent
\cite{cf:carinci1} for a related time dependent problem
with applications to hysteresis. The even $d$ case is instead motivated
by a variety of real world phenomena (e.g., laser instabilities
\cite{cf:pitchfork}, and we suggest to consult the introduction of
\cite{cf:arecchi} for an overview of applications) and, from a more
theoretical viewpoint, by the analysis of anomalous fluctuations at
criticality and here again the analysis
of \textsl{nearly-critical} systems naturally leads to consider
instances like the one in the right inset of Figure~\ref{fig:2}.
Regardless of $d$ being even or odd, and referring to the right insets
of Figure~\ref{fig:1} and
Figure~\ref{fig:2},
for $\mu>0$, respectively $\mu<0$, the fixed point $0$ becomes
linearly stable, respectively unstable,
but observe that the term containing $\mu$ vanishes as $\gep\searrow0$.\looseness=-1
%

%

As we are going to argue further in Section~\ref{sec:gen}, $\mu\neq0$ is
a very relevant generalization, but for the sake of clarity
let us stick to the case $\mu=0$ still for a while. Actually,
the random variable $T_d$ is a limiting universal random variable
behind a very basic mechanism due to
the interaction of noise and non-linearity: it is the first and foremost
quantifier of how a weak stochastic perturbation makes a diffusion go
through a saddle point or how it leads to the escape from a \textsl
{degenerate} unstable point. And in fact there have been several
attempts to determine fine properties of the distribution of $T_d$
in the literature beyond computing the first two moments (see, e.g., \cite{cf:pitchfork,cf:CFB}), but the results appear to be
confined to
uncontrolled approximations and numerical observations that fail, in
particular, to capture for example
the probability of observing large and small values of $T_d$.
The purpose of this work is to present a rigorous treatment of the convergence
statement $\cL-\lim_{\gep\searrow0} \gep^{2(d-2)/d} \tau_{a, \gep
} = T_d$: not surprisingly we will
see that $T_d$ can be directly characterized (recall we are still only
talking about the case $\mu=0$) as the explosion time for
%
\begin{equation}
\label{eq:limSDE} \mathrm{d} Y_t = c_d
Y_t^{d-1} \dd t + \ddd W_t ,
\end{equation}
with $Y_0=0$ for $d$ even and $Y_0=-\infty$ (see Section~\ref{sec:main}
for a precise definition) for $d$ odd ($c_3=1/2$ and $c_4=1$ in the
examples in this introduction).
Moreover, our purpose is to obtain sharp quantitative estimates on the
law of $T_d$, notably
sharp estimates on the probability of observing small and large values
of $T_d$ and regularity estimates on the density of $T_d$.

\subsection{A generalized set-up: Near critical cases}
\label{sec:gen}
As already mentioned above and visually presented in Figures \ref{fig:1} and \ref{fig:2}, it is very natural to consider the
more general set-up
of considering $\gep$-dependent potentials, covering thus
nearly-critical situations.
This actually leads to families of limit transit or escape times
associated to
equations of the form
%
\begin{equation}
\label{eq:limSDEgen} \mathrm{d} Y_t = \sum
_{i=1}^d c_i Y^{i-1} \dd t
+ \ddd W_t .
\end{equation}
While our approach can cover all these cases, the analysis would be
heavy and not particularly transparent.
We therefore decide to deal with the cases in which only the coefficients
$c_d =1$ and $c_2=-\mu$ are non-zero for $d$ even, and $c_d =1/2$ and
$c_1=-\mu/2$ are non-zero for $d$ odd.
In fact, in some cases full details of the analysis will be given only
in the cases $d=3$ and $d=4$. The arising limit
SDEs for which we will be interested in the explosion times, starting
respectively from $-\infty$ and $0$, are
%
\begin{equation}
\label{eq:limSDEtreated} \mathrm{d} Y_t = \frac{1}2
\bigl(Y_t^{2}- \mu \bigr) \dd t + \ddd W_t
\quad \mbox{and}\quad  \mathrm{d} Y_t = \bigl(Y_t^{3}-
\mu Y_t \bigr) \dd t + \ddd W_t .
\end{equation}
Note that the drift part of these two SDEs are the normal form,
respectively, of a saddle-node
bifurcation and of a subcritical pitchfork bifurcation, and of course
$\mu$ is the bifurcation parameter.

It is worth pointing out that generalizing the odd $d$ case to go
toward \eqref{eq:limSDEgen} is really just a matter of
going in detail through the analysis, whereas generalizing the analysis
in the even $d$ case
to non-even limit potentials -- note that the potential in the
right-hand side of \eqref{eq:limSDEtreated}, $-x^4/4+\mu x^2/2$, is
indeed even -- requires a slightly different analysis,
because then the explosion to $\pm\infty$ does not happen with the
same probability (nevertheless, we are able to adapt our
approach also to these cases \cite{cf:GM2}).

Let us conclude the \hyperref[sec:intro]{Introduction} by remarking that in the case $d=2$,
for example $U(x)=-x^2/2$, the point $0$ is hyperbolic unstable.
The noise induced escape from an hyperbolic unstable point happens on
times that to leading order behave like $\log(1/\gep) /\vert
U''(0)\vert$ as can be easily guessed by solving explicitly the
linearized equation. Therefore $T_2$ is, to leading order, just a
positive constant: a detailed treatment can be found for example
in \cite{cf:Bak}, along with the treatment of the subleading
correction, which is random.

\section{Set-up and main results}
\label{sec:main}

\subsection{General set-up and rescaling}
\label{sec:cubic}
For $\mu\in\bbR$ we consider a family of $C^1$ potentials $\{U\me\}
_{\gep>0}$,
with $U\me'(\cdot)$ locally Lipschitz, for which
further assumptions will be given just below (of course the four
potentials in Figures \ref{fig:1} and
 \ref{fig:2} fall into the realm of our analysis).
We consider the strong solution $X$ to the stochastic differential equation
%
\begin{equation}
\label{eq:diffusion} X_0 = x_0^{(\varepsilon)},\qquad  \mathrm{d}
X_t = -(U\me)'(X_t) \dd t + \varepsilon\dd
W_t.
\end{equation}
%

Equation \eqref{eq:diffusion} can be rewritten by performing a rescaling.
For this let us introduce $Y_t := \varepsilon^{-2/d}
X_{\varepsilon^{-2(d-2)/d}t}$,
$B_t = \varepsilon^{(d-2)/d} W_{\varepsilon^{-2(d-2)/d}t},\ t \ge0$
(so $B_\cdot$ is also a standard BM), $y_0^{(\varepsilon)} :=
\varepsilon^{-2/d} x_0^{(\varepsilon)}$, and
$V\me(y) = \varepsilon^{-2} U\me(\varepsilon^{2/d}y) $ so that
$(V\me)'(y) = \varepsilon^{2(1-d)/d} (U\me)'(\varepsilon^{2/d}y)$.
Then $Y$ is the strong solution to the stochastic differential equation
%
\begin{equation}
\label{eq:rescdiff} Y_0 = y_0^{(\varepsilon)}, \qquad \mathrm{d}
Y_t = -(V\me)'(Y_t) \dd t + \ddd
B_t .
\end{equation}

The scaling exponent has been chosen in particular so that $V_{\mu,
\gep}$ has a non-trivial limit as $\gep\searrow0$ and for this we
give the following assumption.

\begin{assume}
\label{assume:V}
Let
%
\begin{equation}
\label{eq:assumeV} V_{\mu}(y) := %
\cases{ -\displaystyle \frac{1}{2d}
{y^d} +\displaystyle \frac{\mu}{2} y &\quad \mbox{if $d$ is odd,}
\vspace*{2pt}\cr
-
\displaystyle \frac{1}d {y^d} +\displaystyle \frac{\mu}{2}y^2 &\quad \mbox{if
$d$ is even.} } %
\end{equation}
We assume that the family of (differentiable, with locally Lipschitz
derivative) functions $\{U\me\}_{\gep>0}$~-- recall that $V\me(\cdot)=\gep^{-2} U\me( \gep^{2/d}\cdot)$
-- is such that, for any $A>0$
%
\begin{equation}
\label{eq:closetoV} \lim_{\gep\searrow0} \llVert V\me- V_{\mu}
\rrVert _{\infty,[-A,A]} = 0 ,
\end{equation}
with, for $B \subset\bbR$, $\Vert f\Vert_{\infty,B} := \sup_{x \in
B} |f(x)|$.
We assume in addition that there exists a $b\in(0, \infty]$, $A>0$
and an increasing continuous function $\psi\dvtx (0, \infty)\to(0, \infty
)$ with the property
$\int_0^\infty\ddd x /\psi(x)<\infty$, such that when $\gep$ is
small enough,
%
\begin{equation}
\label{eq:outoforigin} \psi\bigl(\vert y\vert\bigr) \le %
\cases{ -(V
\me)'(y) &\quad \mbox{if $d$ is odd,}
\cr
-\sign(y) (V
\me)'(y) &\quad \mbox{if $d$ is even,} } %
\end{equation}
for every $y$ such that $ |y| \in[A, b \varepsilon^{-2/d})$.
\end{assume}

Recall that we focus on how fast our initial diffusion $X$
(with potential $U\me$) travels from $x^{(\gep)}_0<0$, notably in the
case in which $ \lim_{\gep\searrow0} x^{(\gep)}_0 <0$,
to $a>0$ (for $d$ odd),
or how long it takes to go from $0$, or very nearby, to $\pm a$ (for
$d$ even). Let us be precise about the initial condition:

\begin{assume}
\label{assume:start2}
For $d$ odd and $b\in(0, \infty)$ (recall that $b\in(0 , \infty]$
is chosen as in Assumption~\ref{assume:V})
we require
%
\begin{equation}
\label{eq:start2odd} \liminf_{\varepsilon\searrow0} x_0^{(\varepsilon)}>
-b \quad \mbox{and}\quad  \gep ^{-2/d}x_0^{(\varepsilon)}=y_0^{(\varepsilon)}
\stackrel{\gep \searrow0}\longrightarrow -\infty
\end{equation}
and if $b= \infty$ it is sufficient to require only the second of the
two conditions in \eqref{eq:start2odd}.
For $d$ even instead we require
%
\begin{equation}
\label{eq:start2even} x_0^{(\varepsilon)}= \mathrm{o} \bigl(
\gep^{2/d} \bigr) .
\end{equation}
\end{assume}

We are after\vspace*{-1pt}
%
\begin{equation}
\label{eq:tau} \tau_{a, \gep}(X) := %
\cases{ \inf \{t\dvt
X_t=a \}  & \quad \mbox{if $d$ is odd,}
\cr
\inf \{t\dvt  \vert
X_t\vert=a \}  &\quad  \mbox{if $d$ is even.} } %
\end{equation}
It is necessary to assume
%
\begin{equation}
\label{eq:b>a} b \ge a ,
\end{equation}
simply because we make no hypothesis on $\{U\me\}(x)$ for $\vert
x\vert>b$.

Equation (\ref{eq:closetoV}) expresses that, as $\gep\searrow0$,
$U\me$ has a very precise limiting behaviour in\vspace*{1pt} any small neighborhood
of $0$ of the form $[-A \varepsilon^{2/d}, A\varepsilon^{2/d}]$.
Equation (\ref{eq:outoforigin}) instead says that $U'\me$ is \textsl
{sufficiently superlinear} with the correct sign in $(-b,b)\setminus
[-A \varepsilon^{2/d}, A\varepsilon^{2/d}]$: this guarantees that, in
the odd $d$ case, the rescaled $Y$ diffusion reaches $-A$ in a finite
time -- recall that the initial condition is \textsl{infinitely far}
from the origin, cf. \eqref{eq:start2odd} -- and that it will escape
\textsl{infinitely far} to the right
once $A$ is reached, again in a finite time. Analogous observations
hold for the even $d$ case.

\subsection{Main results: Convergence and sharp estimates on
generating~functions}
\label{sec:mainres}
Given a random variable $Z\ge0$, we write $\Phi_Z(\gl):=\bbE[ \exp
(\gl Z)]$:
$\Phi_Z(\cdot)$ is the moment generating function of $Z$ (with slight
abuse of notation we use this terminology without further assumptions
on $Z$, notably without
assuming the existence of the moments). Note also that $\Phi_Z(- \cdot
)$ is the Laplace transform of $Z$, but later on we will often use
generating function and Laplace transform as synonymous. Of course
$\Phi_Z(\gl)\le1$ for $\gl\le0$ and $\Phi_Z(\cdot)$ is convex
and non-decreasing. We set
%
\begin{equation}
\label{eq:lambda0} \gl_0=\gl_0(Z):= \sup \bigl\{\gl\dvt
\Phi_Z(\gl)< \infty \bigr\} ,
\end{equation}
so $\Phi_Z(\cdot)$ is well defined and analytic in $\{\gl\in
\bbC\dvt  \Re(\gl) < \gl_0\}$.

We are now ready to state the basic convergence result.

\begin{theorem}
\label{th:lim3}
Under Assumptions \ref{assume:V} and \ref{assume:start2}, we have
%
\begin{equation}
\label{eq:lim} \cL- \lim_{\gep\searrow0} \varepsilon^{2(d-2)/d}
\tau_{a, \gep}(X) =: T_{d, \mu} .
\end{equation}
\end{theorem}

Let us start with the estimates on the moment generating function: of course
Theorem~\ref{th:lim3} is equivalent to
%
\begin{equation}
\label{eq:limPhi} \lim_{\gep\searrow0} \bbE \bigl[\exp \bigl(\gl
\varepsilon^{2(d-2)/d} \tau_a(X) \bigr) \bigr] =
\Phi_{T_{d, \mu
}}(\gl) ,
\end{equation}
for every $\lambda< 0$. But we will see that $\Phi_{T_{d, \mu}}(\gl
)< \infty$ also for some $\gl>0$, that is, with the notation
that we have introduced, $\gl_0(T_{d, \mu})>0$. In order to state our
result on
$\gl_0(T_{d, \mu})$ and on the behavior of $\Phi_{T_{d, \mu}}(\gl
)$ near $\gl_0(T_{d, \mu})$, we
introduce
the Schr\"odinger operator $L$ with domain $C^2(\bbR; \bbR)$:
%
\begin{equation}
\label{eq:Schr} Lu(x) = -\tfrac{1}2 u''(x)
+q_{d, \mu}(x) u(x) ,
\end{equation}
where $q_{d, \mu}(\cdot)$ is a polynomial function, more precisely:
%
\begin{equation}
\label{eq:qdmu} q_{d, \mu}(x) := \tfrac{1}{2} \bigl(
V'_\mu(x) \bigr)^2- \tfrac{1}{2}
V''_\mu(x) .
\end{equation}
In particular, we have
%
\begin{equation}
\label{eq:Schr-Q} q_{3, \mu}(x) = \tfrac{1}2 x+ \tfrac{1}8
\bigl(x^2-\mu \bigr)^2 \quad \mbox{and}\quad  q_{4,
\mu}(x) =
\tfrac{1}2 \bigl(x^3-\mu x \bigr)^2 +
\tfrac{3}2 x^2 .
\end{equation}
Classical deep results, see, for example, \cite{cf:CL,cf:T1} and
Section~\ref{sec:ode}, ensure that the equation $Lu=\gl u$ has a
(classical) solution
that is in $\bbL^2(\bbR; \bbR)$ if and only if
$\gl= \widetilde\gl_j$, with $\widetilde\gl_0< \widetilde\gl_1<
\cdots\,$. Therefore, $\tilde\gl_0$ is the
bottom of the spectrum of $L$. This fact will be used for $d$ odd.

A slightly different result will be needed for even $d$ and this is
connected to the
fact that the question we ask differs according to the parity of $d$.
When $d$ is even, we consider in fact  the spectrum
of $L$ on the domain $[0, \infty)$, or equivalently on $(-\infty,
0]$, with the boundary condition $u(0)=1$ and $u'(0)=0$;
classical results,
see, for example, \cite{cf:CL,cf:T1} and Section~\ref{sec:ode2},
warrant that the spectrum (in $\bbL^2[0, \infty)$) of $L$
is discrete and the eigenvalues $\widetilde\gl_j$ form an increasing
sequence of real numbers like in the previous case.
Actually, in both the even and odd case we will be mainly
interested in $\widetilde\gl_0$ (however note that $\widetilde\gl
_1$ does play a r\^ole in the precision of our approximation in
Proposition~\ref{th:regular} below).
It is practical to introduce the richer notation
$\widetilde\gl_{0,d}, \widetilde\gl_{1,d}$ and the reader should
keep in mind that the even and odd cases correspond  to different
spectral problems,
both because the two Schr\"odinger operator differ and because the domains and boundary conditions also differ.

\begin{theorem}
\label{th:Phiright}
$\gl_0:=\lambda_0(T_{d, \mu})$ is (strictly) positive and it
coincides with $\widetilde\gl_{0,d}$. Moreover, there exists
a positive constant $C_{d, \mu}$
such that
%
\begin{equation}
\label{eq:Phi3.0} \Phi_{T_{d, \mu}}(\gl)\stackrel{\gl\nearrow
\lambda_0}\sim\frac
{C_{d, \mu}}{\lambda_{0} -\gl} .
\end{equation}
Actually $\Phi_{T_{d, \mu}}(\cdot)$ can be extended to the whole of
$\bbC$ as a meromorphic function:
in particular, \eqref{eq:Phi3.0} is therefore saying that $\Phi
_{T_{d, \mu}}(\cdot)$ has a simple pole at $\gl_0$
with residue $-C_{d, \mu}$.
\end{theorem}

See Corollary~\ref{cor:lambda0} and Corollary~\ref{cor:lambda0quart}
for the precise values
of $C_{d, \mu}$, respectively in the odd and even cases: full details
are given in the cases $d=3$ and $4$, but
the generalization is straightforward.

Instead
the asymptotic behavior for $\gl\to-\infty$ of $\Phi_{T_{d, \mu
}}(\gl)$ -- the next result --
is even more explicit, but in this
case the expressions for general $d$ are rather cumbersome (the
expression for the leading term inside the exponential can nonetheless
be generalized in a straightforward way, see below Corollary~\ref
{cor:firstorder} and Remark~\ref{rem:replace} for odd $d$, and
Corollary~\ref{cor:firstorder-quart} for even $d$). So the precise
statement is restricted to $d=3$ and $4$.

\begin{theorem}
\label{th:Phileft}
For $\gl\to-\infty$, we have
%
\begin{equation}
\label{eq:Phi3.1} \Phi_{T_{3, \mu}}(\gl) = \bigl(1+ \mathrm{O} \bigl(\vert
\lambda\vert ^{-1/4} \bigr) \bigr)\exp \bigl(- C_{3/4} \vert\gl
\vert^{3/4}-\mu C_{1/4} \vert\gl\vert^{1/4} \bigr) ,
\end{equation}
where $C_{3/4}= 3 \Gamma (-(3/4) )^2/( 2^{9/4} \sqrt{2 \uppi
})$ and $C_{1/4}=2^{1/4} \Gamma(3/4)^2\sqrt{2/\uppi}$.
Moreover, in the same limit
%
\begin{equation}
\label{eq:Phi3.2} \Phi_{T_{4, \mu}}(\gl) = 2^{-1/4} \vert\gl
\vert^{1/4} \bigl(1+ \mathrm{O} \bigl(\vert\lambda\vert^{-1/3}
\bigr) \bigr) \exp \bigl(-C_{2/3}\vert\gl\vert^{2/3} - \mu
C_{1/3} \vert\gl \vert^{1/3} - \tfrac{1} 6
\mu^2 \bigr) ,
\end{equation}
with $C_{2/3}$ and $C_{1/3}$ positive constants explicitly given in
\eqref{eq:forPl} in terms of elliptic integrals
of first and second kind.
\end{theorem}

\subsection{Tail probabilities, existence and smoothness of density}
\label{sec:tail}

By Tauberian arguments one can extract from Theorem~\ref{th:Phiright}
the behavior of $\bbP(T_{d, \mu}>t)$ for $t$ large
and Theorem~\ref{th:Phileft} yields precise Laplace estimates on $\bbP
(T_{d, \mu}<t)$ for $t \searrow0$. Here is the result:
%

\begin{cor}
\label{th:tails}
We have
%
\begin{equation}
\label{eq:largeT} \lim_{t \to\infty} \frac{1}t \log\bbP
(T_{d, \mu}>t ) = -\lambda_0(T_{d, \mu}) ,
\end{equation}
and
%
\begin{equation}
\label{eq:smallT} \lim_{t \searrow0}{t^{d/(d-2)}}\log\bbP
(T_{d, \mu}<t ) = - a_{d} ,
\end{equation}
with
%
\begin{equation}
a_{d} := \frac{d-2}d \biggl(\frac{d C_{d/(2(d-1))}}{2(d-1)}
\biggr)^{2\vafrac{d-1}{d-2}} .
\end{equation}
\end{cor}

Equation \eqref{eq:largeT} is a direct consequence of Theorem~\ref
{th:Phiright} and the Tauberian result (\cite{cf:nakagawa}, Th. 3).
Instead \eqref{eq:smallT} follows from Theorem~\ref{th:Phileft} and
de Bruijn's Tauberian Theorem (\cite{cf:RegVar}, Th. 4.12.9).

Of course
Corollary~\ref{th:tails} loses quite some information with respect to
Theorem~\ref{th:Phiright} and Theorem~\ref{th:Phileft}.
Theorem~\ref{th:Phiright} actually suggests that the right-tail of
$T_{d,\mu}$ should be close to the tail of
an exponential random variable of parameter $\lambda_0$
translated by
$(1/\lambda_0)\log(C_{d, \mu}/\lambda_0)$. Theorem~\ref
{th:Phileft} yields sharp asymptotic, notably $\mu$ dependent,
behaviors of which
there is no trace in \eqref{eq:smallT}.
%
\begin{figure}

\includegraphics{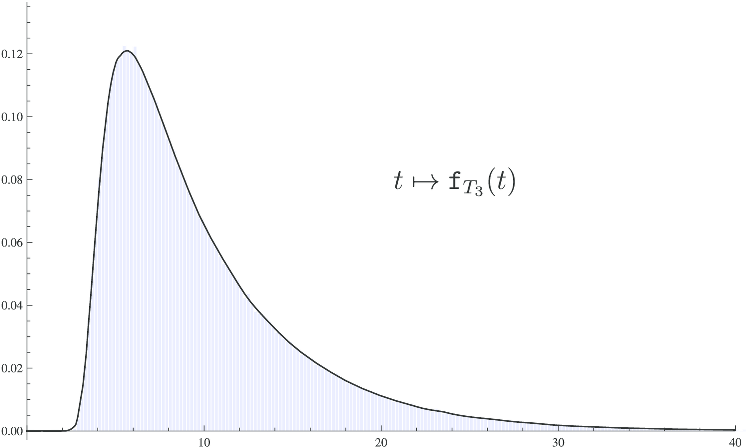}

\caption{The density of $T_{3}$, from the next section it is rather
called $T_{3,0}$, plotted via the histogram representation of a sample
of $5\times10^5$
(independent) realizations of the random variable. The sample has been obtained
by exploiting the representation formula in the first line of \protect\eqref
{eq:Ty0}, but it is particularly challenging to simulate correctly
extreme values (see however the $\log$-plot in the inset of Figure~\protect\ref{fig:compare5}). The empirical mean of the sample we have used is
9.939$\ldots\,$,
against the exact value 9.952$\ldots\,$, and the empirical standard
deviation is 5.70$\ldots$ against
the exact value 5.74$\ldots$ (the first two moments of $T_3$ have
explicit expressions involving $\Gamma$ functions
\cite{cf:SH}).}\label{fig:density}
\end{figure}

\begin{figure}

\includegraphics{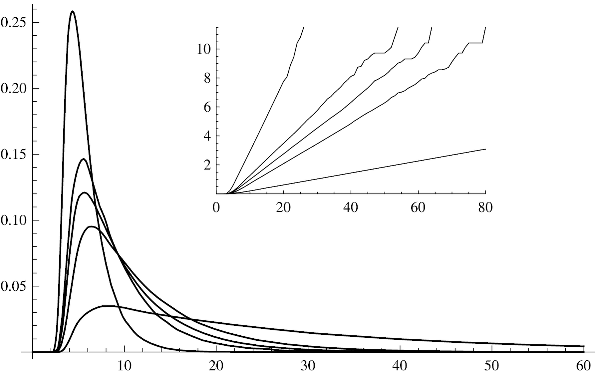}

\caption{The main plot is the one of $t \mapsto\dens_{T_{3,\mu}}$,
for $\mu=-1, 0, -0.2, 0, 0.2, 1$: for $\mu=-1$, respectively for $\mu=+1$,
$T_{3,\mu}$ is more concentrated (resp., more spread out). In the
inset instead there is $t \mapsto-\log\bbP(T_{3,\mu}>t)$, for the
same values of $\mu$. The curves in the main figure have been obtained
by a smoothing procedure on (five) samples each of size $10^5$.
The same samples have been used for the second figure and the plot is
the linear interpolation of the discrete graph obtained for
$t$ that takes integer values.}\label{fig:compare5}
\end{figure}

We can actually strongly improve (\ref{eq:largeT}) of Corollary~\ref{th:tails},
and the
first step 
goes through establishing the existence
of the density $\dens_{T_{d, \mu}}(\cdot)$ of $T_{d, \mu}$ (see Figures~\ref{fig:density}, \ref{fig:compare5}). We can do
better than this, in the sense that we can establish
not only the existence of the density, but also its analyticity
properties and its asymptotic behaviour at $\infty$.
In fact, we will first establish that,
with the standard notation for the characteristic function $\gp_X(s):=
\bbE\exp(\mathrm{i}sX)$ of a random variable $X$, there exists
$c>0$
such that as $s \to\pm\infty$
%
\begin{equation}
\label{eq:forfT} \bigl\vert\gp_{T_{d, \mu}}(s)\bigr\vert= \mathrm{O} \bigl(\exp \bigl(-c
\vert s \vert ^{\afrac{d}{2(d-1)}} \bigr) \bigr) ,
\end{equation}
which we prove in Corollary~\ref{cor:firstorder} and the discussion
following it for $d=3$, and in Corollary~\ref{cor:firstorder-quart}
for $d=4$.

It is a very standard (Fourier analysis) result that \eqref{eq:forfT} entails
that the density $\dens_{T_{d, \mu}}$ exists and it can be chosen to
be $C^\infty$. But in Corollary~\ref{cor:firstorder}
and in Corollary~\ref{cor:firstorder-quart} we actually prove a
result that is substantially stronger than the bound \eqref{eq:forfT},
in the sense that we know the asymptotic behavior of
$\Phi_{T_{d, \mu}}(\cdot)$ along any ray in the complex plane and
therefore we know in which sector the Laplace transform
decays to zero, and how fast, at infinity. This
implies both a stronger result on the regularity of $\dens_{T_{d, \mu
}}$ and, coupled to Theorem~\ref{th:Phiright}, a sharp result
on $\dens_{T_{d, \mu}}(t)$ for $t\to\infty$: this is the content of
the next statement.

\begin{proposition}
\label{th:regular}
$\dens_{T_{d, \mu}}(\cdot)$ is real analytic except at $0$ and it
can be extended to an analytic function
in the cone
%
\begin{equation}
\biggl\{z \in\bbC\dvt  \Re(z)>0, \bigl\vert\mathrm{arg}(z)\bigr\vert< \uppi \biggl(
\frac{1}2- \frac{1}d \biggr) \biggr\} .
\end{equation}
Moreover for $t\to\infty$
%
\begin{equation}
\label{eq:f_T-sharp} \dens_{T_{d, \mu}}(t) = C_{d,\mu} \exp \bigl(-
\lambda_0(T_{d, \mu
}) t \bigr) + \mathrm{O} \bigl( \exp(- b
t) \bigr) ,
\end{equation}
for any choice of $b\in(\tilde\gl_{0,d}, \tilde\gl_{1,d})=
(\lambda_0(T_{d, \mu}), \tilde\gl_{1,d} )$.
\end{proposition}

The substantial difference between obtaining $t\to\infty$
and $t \searrow0$ estimates on the density is that
in the first case the leading behavior is directly linked to a pole of
the Laplace transform of the density
(see \cite{cf:FGD} for more details), while in the second case an
essential singularity enters the game. As a matter of fact
it is not difficult to show that $\dens_{T_{d, \mu}}(t)$ goes quickly
to zero as $t \searrow0$, but even simply extending
\eqref{eq:smallT} to $\dens_{T_{d, \mu}}(t)$ appears to be rather
challenging: this can be approached by proving a \textsl{slow-decrease}
property for $\dens_{T_{d, \mu}}(\cdot)$ near zero, so that the
Tauberian Theorem (\cite{cf:RegVar}, Th. 4.12.11), would apply, but we do
not have such an estimate.

%

\subsection{Organization of the paper, with a sketch of the approach}
In Section~\ref{sec:mart}, we apply results from martingale theory to
prove Theorem~\ref{th:lim3} along with
a representation of the Laplace transform of the law of $T_{d, \mu}$,
that is $\Phi_{T_{d, \mu}}(\gl)$ for $\gl\le0$,
as ratio of asymptotic values of the solution $f_\gl(\cdot)$ of a
suitable ODE, for $\gl\le0$. The representation formula has the form
%
\begin{equation}
\label{eq:reprfor} \ee \bigl[\exp(\lambda T_{d,\mu}) \bigr] =
\frac{f_{\lambda
}(y_0)}{f_{\lambda}(+\infty)} ,
\end{equation}
with $y_0=-\infty$ for odd $d$ and $y_0=0$ for even $d$.

In Section~\ref{sec:ode}, we make a thorough analysis of the ODE
solved by $f_\gl(\cdot)$, for odd $d$
and for the sake of clarity we give full details only for $d=3$. Via a
standard transformation such ODE
is mapped into a Schr\"odinger equation for which the analysis is
carried out. It is at this level that the spectral properties of the arising
Schr\"odinger operator play a r\^ole, but the questions we are asking
are rather different from the standard ones that are typically
addressed for such an operator. In fact, much of the analysis takes
place out of the spectrum, so for functions that are not in $\bbL^2$.
Moreover, we analyze the solutions of the arising Schr\"odinger
equation to establish that the right-hand side of
\eqref{eq:reprfor} is the ratio of entire functions (of $\gl$), and
it is therefore meromorphic, so that, in particular, \eqref
{eq:reprfor} holds also for $\gl\in\bbC$ with $\Re(\gl)$ smaller
than the first
pole of the function in the right-hand side of
\eqref{eq:reprfor} (and such a pole is real: actually, the set of the
poles of the right-hand side of \eqref{eq:reprfor} coincides with the spectrum
of the Schr\"odinger operator we are working with). By exploiting
\eqref{eq:reprfor}, extended as explained to complex values of $\gl$,
we will then obtain
the proofs of Theorem~\ref{th:Phiright}, of
Theorem~\ref{th:Phileft} and of Proposition~\ref{th:regular}, for $d$ odd.

In Section~\ref{sec:ode2}, we go again through the arguments for the
case of even $d$. 

\section{Martingales, Laplace transforms, and convergence}

\subsection{Martingales and Laplace transforms}
\label{sec:mart}

One should note that most objects of interest introduced in the
paragraph below such as functions $u , s_{\gep}$ or processes $Y,M$ in
fact depend on $V\me$ (in particular they all depend on the two
parameters $\mu$ and $\gep$). For bookkeeping purposes this will not
appear in our notation but the reader should keep this dependency in mind.

For $t \ge0$, write $\cF_t = \sigma(B_s, s \le t)$. By (\ref
{eq:rescdiff}) and a direct application of It\^o's formula, if the
function $u \dvtx [0, \infty) \times\R\to\R$ is $C^2$ and satisfies
%
\begin{equation}
\label{eq:edpuY} \partial_t u (t,y) - \partial_y u(t,y) V
\me'(y) + \tfrac{1}{2}\, \partial^2_y
u (t,y) =0,
\end{equation}
then
%
\begin{equation}
\label{eq:locmartY} u(t,Y_t) = u \bigl(0,y_0^{(\varepsilon
)}
\bigr) + \int_0^t u'(s,Y_s)
\,\mathrm{d}B_s.
\end{equation}
Moreover, if $\tau$ is any $(\cF_t)$-stopping time then $ ( M_t
: = u(t \wedge\tau, Y_{t \wedge\tau}) )_{t \ge0}$ is a $(\cF
_t)$-local martingale, whose quadratic variation at $t \ge0$ equals
$ \int_0^{t \wedge\tau} (u'(s,Y_s))^2 \dd s$. Here are some standard
facts (for a proof see, e.g., Theorem~4.7 in \cite{cf:LG}):

\begin{proposition}
\label{cl:truemart}
If, for all $t \ge0$
%
\begin{equation}
\ee \biggl[ \int_0^{t \wedge\tau} \bigl(u'(s,Y_s)
\bigr)^2 \,\mathrm{d}s \biggr] < \infty,
\end{equation}
then $M$ is a square integrable martingale.
Moreover, if
%
\begin{equation}
\ee \biggl[ \int_0^{\tau} \bigl(u'(s,Y_s)
\bigr)^2 \,\mathrm{d}s \biggr] < \infty,
\end{equation}
then $M$ is bounded in $\mathbb{L}^2$, hence uniformly integrable.
\end{proposition}

We recall also the Doob's Optional Stopping Theorem: if $M$ is a
uniform integrable martingale
and $\tau$ a stopping time, then $M_\tau\in\bbL^1$ and
$\bbE[M_\tau]=\bbE[M_0]$. 

\subsubsection*{A first application of Proposition \protect\ref{cl:truemart}:
A localization lemma}

Proposition~\ref{cl:truemart} directly yields exit probabilities from
a strip. This is useful only
when $d$ is odd: in the even $d$ case
we know a priori that $X_{t\wedge\tau_{a, \gep}}$ stays in $[-a,a]$,
but in the odd $d$ case
$X_{t\wedge\tau_{a, \gep}}$ can in principle go arbitrarily far off
to the left, even where we cannot control it anymore (i.e., to the left
of $-b$). So, let us fix $d$ odd and argue that this happens with
negligible probability: we will first do this
under the assumption $y_0^{(\varepsilon)}\stackrel{\gep\searrow
0}{\longrightarrow} y_0\in(-\infty, 0)$. Note that this is in
contrast with Assumption~\ref{assume:start2}, that requires
$y_0=-\infty$: we will in fact first treat the case $y_0$ finite
and then show how to let $y_0 \to-\infty$.
\begin{figure}

\includegraphics{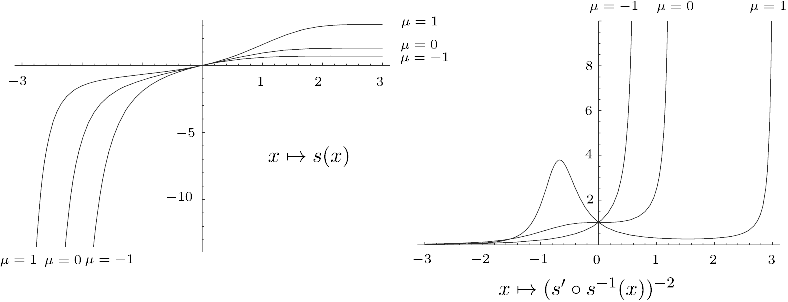}

\caption{The plot of the scale function $s(\cdot)$ ($d=3$)
and the function $(s'\circ s^{-1}(\cdot))^{-2}$, both for $\mu=-1$,
$0$ and $+1$. The latter function is the one appearing in \protect\eqref{eq:Ty0}:
while for $\mu<0$ such a function is increasing and the large values
are close to the divergence, that is close to $x=s(\infty)$, for $\mu>0$
another bump appears and this accounts for the presence of a valley in
the potential that (weakly) traps the diffusion and yields longer
transit times.}\label{fig:s}
\end{figure}

In order to control the probability of an excursion far off to
the left, we introduce a new stopping time for $d$ odd:
for $ -b \gep^{-2/d} \le l_{\gep}<y_0^{(\varepsilon)}<r_{\gep} \le
b \gep
^{-2/d}$ let us set
%
\begin{equation}
\label{eq:sigma} \sigma= \sigma_{l_{\gep},r_{\gep}}(Y) := \inf \bigl\{t \ge0 \dvt
Y_t \in \{l_{\gep}, r_{\gep}\} \bigr\} .
\end{equation}
It is easily checked that $\gs$ is a stopping time and that it is
(exponentially) integrable.
Let us introduce also (see Figure~\ref{fig:s})
%
\begin{equation}
\label{eq:fmart} s\e(y) := \int_0^y \exp
\bigl(2 V\me(u) \bigr)\,\mathrm{d}u,\qquad  y \in\R.
\end{equation}
One easily checks that $s\e$ satisfies (\ref{eq:edpuY}), thus $(s\e
(Y_{t\wedge\sigma}))_{t \ge0}$ is a local martingale.
In addition, for any $t \le\sigma$, $|Y_t| \le\max(|l_{\gep
}|,|r_{\gep}|)$ so that
$|s\e'(Y_t)| \le C$. Thus
%
\begin{equation}
\label{eq:fmartc} \ee \biggl[ \int_0^{\sigma} \bigl(s
\e'(Y_t) \bigr)^2 \dd t \biggr] \le
C^2 \ee [\sigma],
\end{equation}
hence $(s\e(Y_{t\wedge\sigma}))_{t \ge0}$ is a martingale, bounded
in $\mathbb{L}^2$, so
Doob's Optional Stopping Theorem implies $\ee[s\e(Y_{\sigma})] = s\e
(y_0)$. But $Y_\gs\in\{l_{\gep},r_{\gep}\}$ so
%
\begin{equation}
\label{eq:p+} \p(Y_{\sigma} = r\e) = \frac{s\e(y_0^{(\varepsilon)})-s\e
(l_{\gep})}{s\e(
r_{\gep})-s\e(l_{\gep})} =
\frac{\int_{l_{\gep
}}^{y_0^{(\varepsilon)}} \exp
 (2 V\me(u) ) \,\mathrm{d}u}{\int_{l_{\gep}}^{r_{\gep}} \exp
(2 V\me(u) ) \,\mathrm{d}u}.
\end{equation}

From now on, we set
%
\begin{equation}
\label{eq:lere} l\e:= - \gep^{-2/d} \times %
\cases{
(b-x_0)/2 &\quad  \mbox{if }$b< \infty$,
\cr
2x_0 & \quad \mbox{if }$b=\infty$,} \qquad
r\e:= a \gep^{-2/d}.
\end{equation}
We can apply \eqref{eq:p+} to the above set-up,
with $y_0^{(\varepsilon)}= \gep^{-2/d}x_0$, $x_0\in(-b, 0)$,
and $\sigma\e:= \sigma_{l\e,r\e}(Y)$, and find
%
\begin{equation}
\p(Y_{\sigma\e} = r\e) = \frac{\int_{l\e}^{y_0^{(\varepsilon)}}
\exp
(2 V\me(u) ) \dd u}{\int_{l\e}^{r\e} \exp (2 V\me
(u) ) \dd u} \stackrel{\varepsilon\searrow0}
{\longrightarrow} 1 .
\end{equation}
%
Therefore, we have the following lemma.

\begin{lemma}\label{th:local}
For $d$ odd and $l_\gep$ as in \eqref{eq:lere}, we have
%
\begin{equation}
\label{eq:sigmatauY} \lim_{\gep\searrow0} \bbP \bigl( \sigma_{l\e,r\e}(Y)
\ne\tau_{a\varepsilon^{-2/d}}(Y) \bigr) = 0 .
\end{equation}
\end{lemma}

Thus, the convergence in distribution of $\varepsilon^{2(d-2)/d} \tau
_a(X)$ is equivalent to that of  $\varepsilon^{2(d-2)/d}\times\allowbreak \sigma_{-c,a}(X)$,
for a wise choice of the positive constant $c$.


\subsubsection*{Informal martingale approach to the Laplace transform
of $\tau$}
\label{sec:Lapmart}

A classical approach to computing the Laplace transform (or
equivalently the moment generating function) of the distribution of a
hitting time is through martingales.
When interested in the Laplace transform of $\tau_{a \varepsilon
^{-2/d}}(Y) = \varepsilon^{2(d-2)/d} \tau_a(X)$
the game is to find a martingale of the form
%
\begin{equation}
\label{eq:martLap} M_t := \exp(\lambda t) f_{\lambda}(Y_t)
,
\end{equation}
with the additional condition that $(M_{t \wedge\tau_{a \varepsilon
^{-2/d}}(Y)})_{t \ge0}$ is uniformly integrable.

By It\^o's formula, if $f_{\lambda, \gep}$ satisfies
%
\begin{equation}
\label{eq:edoeps} \tfrac{1}{2} f_{\lambda, \gep}''(x)-
(V\me)'(x) f_{\lambda, \gep
}'(x) +\lambda
f_{\lambda, \gep}(x) =0,
\end{equation}
then
%
\begin{equation}
\label{eq:martLapbis} \mathrm{d} M_t = \exp(\lambda t) f_{\lambda, \gep}'(Y_t)
\dd B_t .
\end{equation}
At this stage the analysis of odd and even $d$ differ. Let us consider
the odd $d$ case:
even
for $\lambda< 0$, which we assume, it is in general false that
%
\begin{equation}
\label{eq:Laplacetaueps} \ee \biggl[\int_0^{\tau_{a\gep^{-2/d}}(Y)} \exp(2
\lambda s) \bigl(f_{\lambda, \gep}'(Y_s)
\bigr)^2 \,\mathrm{d}s \biggr] < \infty.
\end{equation}
As a matter of fact \eqref{eq:edoeps} has a two dimensional space of
solutions and it is conceivable that
a choice has to be made at this stage so that \eqref{eq:Laplacetaueps}
holds true, for example that
$\sup_{ y \le a \gep^{-2/d}}\vert f_{\lambda, \gep}'(y) \vert<
\infty$. Let us proceed assuming \eqref{eq:Laplacetaueps}:
we are then dealing with a martingale bounded in $\mathbb{L}^2$, so
(cf. Proposition~\ref{cl:truemart})
$\ee [M_{\tau_{a\gep^{-2/d}}(Y)} ] = M_{0}$. We therefore find
%
\begin{equation}
\label{eq:Laplacetauepsbis} \ee \bigl[\exp \bigl(\lambda\tau_{a\gep^{-2/d}}(Y) \bigr) \bigr]
= \frac{f_{\lambda,
\gep}(y_0^{(\varepsilon)})}{f_{\lambda, \gep}(a \varepsilon
^{-2/d})} .
\end{equation}
Of course, the difficulty is that the ordinary differential equation
(\ref{eq:edoeps}) is not explicit (and its solution $f_{\lambda, \gep
}$ even less so). Indeed the equation depends on $V\me(\cdot)$ which
itself may depend on our parameter $\varepsilon$ in a non-trivial
manner. Thus, we are quite far from a satisfactory formula, or even
from a formula \textsl{tout court}.

Heuristically nonetheless one can formally let $\gep\searrow0$ in
\eqref{eq:edoeps} and in \eqref{eq:Laplacetauepsbis}
(recall $V\me(\cdot) \to V_{\mu}(\cdot)$, and $-2V_{\mu}'(x) =
x^{d-1}-\mu$)
%
\begin{equation}
\label{eq:edolimit} f_{\lambda}''(x)+
\bigl(x^{d-1}-\mu \bigr)f_{\lambda
}'(x) +2\lambda
f_{\lambda}(x) = 0 ,
\end{equation}
which has to be supplied by appropriate boundary conditions, and, if we
call $T_{y_0}$ the limit variable, we should have
%
\begin{equation}
\label{eq:LaplaceTy0} \ee \bigl[\exp(\lambda T_{y_0}) \bigr] =
\frac{f_{\lambda
}(y_0)}{f_{\lambda}(+\infty)} ,
\end{equation}
which is the formal limit of
\eqref{eq:Laplacetauepsbis}. Note that, in view of the right-hand side
of \eqref{eq:LaplaceTy0}, it is sufficient
to determine $f_{\lambda}(\cdot)$ up to a multiplicative constant (we
still have one degree of freedom though!).

Finally, one should not forget that in the odd $d$ case we are really
interested in sending $y_0$ to $-\infty$.

Since making rigorous all the steps we have just outlined does not seem
to be easy, we take
the following alternative path:
\begin{itemize}
\item Instead of working directly with $Y$, we go back to \eqref
{eq:fmart} and work with the martingale $s\e(Y)$, which in
turn can be transformed into
a time-changed Brownian motion, thanks to Dubins--Schwarz Theorem. We do
this for $Y_0=y_0$, or in the slightly generalized case
of $y_0^{(\gep)}$ converging to $y_0$. This gives an amenable formula
for $\tau_{a\gep^{-2/d}}(Y)$
(rather, it gives an amenable formula for
$\gs_{l_\gep, r_{\gep}}(Y)$ in the odd $d$ case, but recall
Lemma~\ref{th:local}). This step is performed in Section~\ref{sec:DS}.
\item We can pass to the limit in this formula, see Section~\ref
{sec:fixedy0}, establishing thus that $\tau_{a\gep^{-2/d}}(Y)$
converges in
law as $\gep\searrow0$ to a limit variable that we call $T_{y_0}$,
and this for every allowed choice of $V\me$.
\item Convergence in law is actually equivalent to the convergence of
$\bbE[ \exp( \gl\tau_{a\gep^{-2/d}}(Y)) ]$ to $\bbE[\exp(\gl T_{y_0})]$
for every $\gl<0$.
But we can compute $\bbE[ \exp( \gl\tau_{a\gep^{-2/d}}(Y)) ]$ by making
a judicious choice of $V\me(\cdot)$, that is simply $V\me(\cdot
)=V_\mu(\cdot)$ (cf. \eqref{eq:assumeV}), so that
\eqref{eq:edoeps} becomes \eqref{eq:edolimit} and we have gotten rid
of the $\gep$ dependence in \eqref{eq:edoeps}.
There is still a priori an obstacle in making the steps \eqref
{eq:martLap}--\eqref{eq:Laplacetauepsbis} rigorous: selecting
the right solution of \eqref{eq:edolimit}, since there is one degree
of freedom. But the crucial condition
\eqref{eq:Laplacetaueps} does require some boundedness condition on
$f'(y)$ and since $f(\cdot)$ is smooth this amounts
to require this for $y\to-\infty$. As a matter of fact solutions to
\eqref{eq:edolimit} can have only certain asymptotic behaviors
(this is one of the instances in which the WKB analysis plays a role),
so actually requiring that
$f(y)$ is bounded as $y \to-\infty$, implies that $f'(y) \to0$, and
in the end we will consider the solution of \eqref{eq:edolimit}
such that $f'(-\infty)=0$ and (say) $f(0)=1$.
\item Finally, in Section~\ref{sec:toinfty}, we will send $y_0\to-\infty$
in the left-hand side of \eqref{eq:Laplacetauepsbis} by performing a
direct SDE estimate
(the non-Lipschitz character of $V_\mu(\cdot)$ makes the time spent
by the diffusion to go from
$-\infty$ to $y_0$ arbitrarily small as $y_0\to-\infty$).
\end{itemize}



\subsection{Scale function, time-changed Brownian motion}
\label{sec:DS}
Let $Z_t = s\e(Y_t)$ with $s_\gep(\cdot)$ given in \eqref{eq:fmart}
and $V(\cdot)=V\me(\cdot)$.
We also write $s(y):=s_0(y)=\int_0^y \exp(2 V_\mu(u)) \dd u$.
Of course $s\e$ is $C^2$ and increasing, and
by (\ref{eq:closetoV}), one directly verifies that, for any ${\mathtt
K} \subset\R$ compact,
$\Vert s\e- s \Vert_{\infty,{\mathtt K}} $ vanishes as $\gep
\searrow0$.
Moreover, for any $a$ such that $0<a\le b$, and for any $y \in
[0,a\gep^{-2/d}]$,
%
\begin{equation}
\label{eq:closetosbis} \bigl\vert s\e(y) - s(y)\bigr\vert\le A \llVert s\e- s\rrVert
_{\infty, [0,A]} + \int_A^{a \gep^{-2/d}} \bigl(\exp \bigl(
2 V\me(u) \bigr) + \exp \bigl( 2 V_\mu(u) \bigr) \bigr) \dd u .
\end{equation}
Thanks to Assumption~\ref{assume:V}, as $A \to\infty$, the
second term in the right-hand side tends to $0$ uniformly in $\gep\in
[0, \gep_0]$ for small enough $\gep_0$. Thus,
$\llVert  s\e-s\rrVert _{\infty, [0,a\gep^{-2/d}]}$ vanishes as $\gep\searrow0$.
A very similar estimate shows
that $\llVert  s\e'-s'\rrVert _{\infty, [-A,a\gep^{-2/d}]}$ vanishes in the same
limit. This actually directly implies for $d$ even
%
\begin{equation}
\label{eq:cts-even} \lim_{\gep\searrow0} \llVert s\e-s\rrVert
_{\infty, [-a\gep^{-2/d},a\gep^{-2/d}]} = 0,\qquad  \lim_{\gep\searrow0} \bigl\llVert s
\e'-s'\bigr\rrVert _{\infty, [-a\gep^{-2/d},a\gep
^{-2/d}]} = 0 ,
\end{equation}
while for $d$ odd \eqref{eq:cts-even} holds with $[-a\gep
^{-2/d},a\gep^{-2/d}]$ replaced
by $[-A,a\gep^{-2/d}]$. The same argument
also implies that
%
\begin{equation}
\label{eq:selim} s\e \bigl(a \gep^{-2/d} \bigr) \stackrel{\gep\searrow0} {
\longrightarrow} \int_{0}^{\infty} \exp
\bigl(2V_\mu( u) \bigr) \dd u = : s(\infty) ,
\end{equation}
and for $d$ even the same holds also for $s\e(-a \gep^{-2/d})$.
On the other hand, by using once again (\ref{eq:outoforigin}), we find
for $d$ odd that for any $l \in(-b,0)$,
%
\begin{equation}
\label{eq:selim2} s\e \bigl(l \gep^{-2/d} \bigr) \stackrel{\gep\searrow0} {
\longrightarrow} \int_{0}^{-\infty} \exp \biggl( -
\frac{u^d}{d} + \mu u \biggr) \,\mathrm{d}u = -\infty.
\end{equation}
%

Recall now (cf. \eqref{eq:fmart}) 
that
$(Z_t)_{t \ge0}$ is a continuous local martingale, started at $s\e
(y_0^{(\varepsilon)})$, and with quadratic variation
%
\begin{equation}
\label{eq:crochetZ} \langle Z \rangle_t = \int_0^t
\exp \bigl(4 V_{\mu}(Y_s) \bigr) \dd s = \int
_0^t \bigl((s\e)' \circ(s
\e)^{-1} \bigr)^2 (Z_s) \dd s ,
\end{equation}
for every $t\ge0$.
Dubins--Schwarz Theorem (see, e.g., Theorem~5.5 in \cite{cf:LG}) then
directly leads to the following lemma.

\begin{lemma} \label{lem:BM}
On an enlarged probability space there exists a standard Brownian
motion $\beta^{(\varepsilon)}$ such that
%
\begin{equation}
\label{eq:Z} Z_t = Z_0 + \beta^{(\varepsilon)}_{\gamma_t^{(\varepsilon)}}
= s\e \bigl(y_0^{(\varepsilon)} \bigr) + \beta^{(\varepsilon)}_{\gamma_t^{(\varepsilon)}}
,
\end{equation}
with
%
\begin{equation}
\label{eq:gamma} \gamma_t^{(\varepsilon)}= \int_0^t
\bigl(g\e(Z_s) \bigr)^{2} \dd s , \quad \mbox{and}\quad  g\e(x) =
(s\e)' \circ(s\e)^{-1} (x) .
\end{equation}
Moreover, the (continuous, increasing) inverse of $\gamma
^{(\varepsilon)}$ is
%
\begin{equation}
\label{eq:A} A_t^{(\varepsilon)}= \inf \bigl\{u \ge0 \dvt
\gamma_u^{(\varepsilon)}>t \bigr\} = \int_0^t
g\e \bigl(\beta^{(\varepsilon)} _s + s\e(y_0)
\bigr)^{-2} \dd s .
\end{equation}
\end{lemma}

\subsection{The limit for a restricted class of initial
conditions} \label{sec:limity0}
\label{sec:fixedy0}
We start with the case when $Y_0$ is a point $y_0\in\bbR$ not
depending on $\gep$.

\begin{assume} \label{assume:start}
Suppose $y_0^{(\varepsilon)}= y_0$ for every $\gep>0$.
\end{assume}


\begin{proposition} \label{prop:limy_0}
Under
Assumption~\ref{assume:start} and assuming without loss of generality
that $y_0\ge0$ in the even case, as $\gep\searrow0$,
$\varepsilon^{2(d-2)/d} \tau_a(X) $ converges in distribution towards
%
\begin{eqnarray}
\label{eq:Ty0} %
T_{y_0} &:=& \int_0^{\tau_{s(\infty)-s(y_0)}(\beta)}
\bigl(s' \circ s^{-1} \bigl(\beta_s+s(y_0)
\bigr) \bigr)^{-2} \,\mathrm{d}s\nonumber
\\[-8pt]\\[-8pt]
& =& \int_{-l}^{s(\infty)} \bigl(s' \circ
s^{-1}(y) \bigr)^{-2} \ell_{\tau
_{(s(\infty)-s(y_0))(\beta)}}^{y-s(y_0)}(
\beta) \dd y ,\nonumber %
\end{eqnarray}
where $l=\infty$ (resp., $l=s(\infty)$) in the odd
(resp., even) case, $\ell_t^x(\beta)$ is the local time of
$\beta$ at level $x$, up to time $t$.
\end{proposition}


\begin{pf}
Let us focus on the odd $d$ case: the even case is
almost identical.
Thanks to (\ref{eq:sigmatauY}), it is equivalent to look at the
convergence in distribution, as $\gep\searrow0$, of $\sigma_{l\e
,r\e}(Y)$. Note that this step is superfluous in the
even $d$ case.
Using that $s\e$ is increasing, and Lemma~\ref{lem:BM}, we find
%
\begin{equation}
\sigma_{l\e, r\e}(Y) \stackrel{\mathrm{(law)}} {=} A^{(\varepsilon)}_{\sigma_{s\e(l\e)- s\e(y_0), s\e(r\e)- s\e
(y_0)}(\beta^{(\varepsilon)})}.
\end{equation}
Since we are only interested in the distribution of the above quantity,
we may as well use a generic Brownian motion $\beta$ for any $\gep$, thus
%
\begin{eqnarray}
\sigma_{l\e, r\e}(Y) &\stackrel{\mathrm{(law)}} {=}& \int
_0^{\sigma_{s\e(l\e)-s\e(y_0),s\e(r\e)-s\e(y_0)}(\beta)} g\e \bigl(\beta_s+s
\e(y_0) \bigr)^{-2} \dd s\nonumber
\\[-8pt]\\[-8pt]
& =& \int_{s\e(l\e)- s\e(y_0)}^{s\e(r\e)-s\e(y_0)} g\e \bigl(x+s\e
(y_0) \bigr)^{-2} \ell_{\sigma_{s\e(l\e)- s\e(y_0),s\e(r\e)-s\e
(y_0)}(\beta)}^x(
\beta) \dd x .\nonumber %
\end{eqnarray}
As $\gep\searrow0$, by (\ref{eq:selim}), $s\e(r\e) \to
s(\infty)$, and by (\ref{eq:selim2}), $ s\e(l\e) \to-\infty$.
Proceeding as for (\ref{eq:cts-even}), we find that for any $l \in
(-b,0)$ and for any $A < s(\infty)$,
%
\begin{equation}
\bigl\llVert (g\e)^{-2} - g^{-2}\bigr\rrVert
_{\infty, [s\e(-\gep^{-2/d}l),A]} \stackrel {\gep \searrow0} {\longrightarrow} 0, \qquad \mbox{with }
g(x) := \bigl(s'\circ s^{-1} \bigr) (x) ,
\end{equation}
see Figure~\ref{fig:s}.
Now, the (deterministic) convergence of $s\e(r\e)$ and $s\e(l\e) $
easily implies that for almost every trajectory of $\gb$,
%
\begin{equation}
\sigma_{s\e(l\e)- s\e(y_0), s\e(r\e)-s\e(y_0)} (\beta)\stackrel {\gep\searrow0} {\longrightarrow}
\sigma_{-\infty,s(\infty
)-s(y_0)}(\beta) = \tau_{s(\infty)-s(y_0)}(\beta) .
\end{equation}
Moreover, $(\ell_t^x(\beta))_{t \ge0, x \in\R}$ is a.s. jointly
continuous in $(t,x) \in\R_+ \times\R$ (see, e.g., Theorem VI.1.7 in
\cite{cf:RY}).
Thus, for fixed positive $K,T$, the above implies that uniformly over
$x \in[-K,K]$, almost surely
%
\begin{equation}
\ell_{\sigma_{s\e(l\e)-s\e(y_0),s\e(r\e)-s\e(y_0)}(\beta)
\wedge T}^{x}(\beta) \stackrel{\gep\searrow0} {
\longrightarrow} \ell_{\tau_{s(\infty)-s(y_0)}(\beta) \wedge T}^x(\beta) ,
\end{equation}
where we have used that a continuous function on an open domain (in
this case of $\bbR^2$) is uniformly continuous
on a compact subset of the domain.
It only remains to note that
%
\begin{equation}
\lim_{T \to\infty} \p \bigl(\tau_{s(\infty)-s(y_0)}(\beta) \ge T/2
\bigr) = 0 \quad \mbox{and}\quad  \lim_{K \to\infty} \p \bigl(\sup \bigl\{|x| \dvt
\ell_T^x > 0 \bigr\} > K \bigr) = 0 ,
\end{equation}
to ensure that, uniformly over $x \in\R$, we have the convergence in
probability
%
\begin{equation}
\ell_{\sigma_{s\e(g\e)-s\e(y_0), s\e(d\e)-s\e(y_0)}(\beta
)}^x(\beta) \stackrel{\varepsilon\searrow0} {
\longrightarrow} \ell _{\tau_{s(\infty)-s(y_0)}(\beta)}^x(\beta) ,
\end{equation}
and therefore also the convergence in probability
%
\begin{equation}
\int_{s\e(g\e)}^{s\e(d\e)} g\e(y)^{-2}
\ell_{\tau_{s\e(d\e
)-s\e(y_0)}(\beta)}^{y-s\e(y_0)}(\beta) \dd y \stackrel{\varepsilon\searrow0} {
\longrightarrow} \int_{-\infty}^{s(\infty)} g(y)^{-2}
\ell_{\tau_{s(\infty
)-s(y_0)}(\beta)}^x(\beta) \dd y .
\end{equation}
This concludes the proof of Proposition~\ref{prop:limy_0}.
\end{pf}

\begin{cor}
\label{cor:edo}
Fix $\lambda\le0$.

For $d$ odd, assume that there exists a solution to
%
\begin{equation}
\label{eq:edofgl} f_{\lambda}''(y) +
2V'_\mu(y)f_{\lambda}'(y)+ 2\lambda
f_{\lambda
}(y) = 0 ,
\end{equation}
such that $\lim_{y \to-\infty} f'_\gl(y)=0$.

For $d$ even instead we consider a solution of \eqref{eq:edofgl} such that
$f'_\gl(0)=0$.

Then in both cases, $\lim_{y \to\infty}f_\gl(y)=f_\gl(\infty)\in
(0, \infty)$ exists and we have
%
\begin{equation}
\label{eq:coredo} \ee \bigl[\exp(\lambda T_{y_0}) \bigr] =
\frac{f_{\lambda}(y_0)}{f_{\lambda
}(\infty)} .
\end{equation}
\end{cor}

Note that, in the odd case, \eqref{eq:coredo} implies the uniqueness
of $f_\gl(\cdot)$ up to a multiplicative
constant (which will be actually chosen in the proofs, by prescribing
the precise asymptotic behavior
of $f_\gl(\cdot)$ near $-\infty$): the existence will be proven in
Section~\ref{sec:ode}.

When $d$ is even, note that $f_\gl(\cdot)$ is an even function, so
$T_{y_0} $ coincides in law with $T_{-y_0}$. Moreover,
we have $T_{0}=T_{d,\mu}$ with the notation in Theorem~\ref{th:lim3}.

\begin{pf*}{Proof of Corollary~\ref{cor:edo}}
Recall our reasoning from paragraph \ref{sec:mart}. We are now in
position to get around the difficulty of the fact that the ODEs (\ref
{eq:edoeps}) may depend on $\gep$ in a non-trivial way.
Indeed from Proposition~\ref{prop:limy_0}, we know that the limit
$T_{y_0}$ \emph{does not depend} on our choice of the family of
potentials $\{V\me\}_{\gep>0}$ (as long as they satisfy Assumption~\ref{assume:V}).
To characterize the distribution of the limit, we may therefore as well choose
$V\me= V_{\mu}$ for any $\gep>0$.
For this particular choice of the family of potentials, the ODE becomes
(\ref{eq:edofgl}). In such a case, for any $\gep>0$
%
\begin{equation}
\ee \bigl[\exp \bigl(\lambda\tau_{a\gep^{-2/d}}(Y) \bigr) \bigr] =
\frac{f_{\lambda
}(y_0)}{f_{\lambda}(a \varepsilon^{-2/d})} .
\end{equation}
Since $\lambda\le0$ we can let $\gep\searrow0$ in the left-hand
side, and by Proposition~\ref{prop:limy_0}, the limit is $\ee[\exp
(\lambda T_{y_0})]$.
This implies that $f_{\lambda}(\cdot)$ has a non-zero limit at
infinity, and that
\eqref{eq:coredo}
holds.
\end{pf*}


\begin{remark} \label{rk:y_0ep}
Results of Proposition~\ref{prop:limy_0}, Corollary~\ref{cor:edo} are
unchanged
when Assumption~\ref{assume:start} is relaxed into the weaker
$y_0^{(\varepsilon)}\to y_0 $ as $\gep\searrow0$.
Indeed, in the above proofs, when before letting $\gep\searrow0$, one
simply need to replace occurrences of $y_0$ with $y_0^{(\varepsilon)}$.
\end{remark}

\begin{remark} \label{rk:McGill}
Even if this meant changing the scale of $\beta$, we could also
have used a result of MacGill \cite{cf:McGill} (important steps in
Macgill's paper are outlined in \cite{cf:RY}, Exercise XI.2.7) which
expresses directly the Laplace transform of an additive functional of
local times, such as the one in
\eqref{eq:Ty0}. The reader who would care to check this alternate
argument would then rather find that
%
\begin{equation}
\ee \bigl[\exp(\lambda T_{y_0}) \bigr] = \frac
{1}{h_{\lambda
}(1)},
\end{equation}
where $h_{\lambda}$ is the unique solution to
%
\begin{eqnarray}
\label{eq:h} &h_{\lambda}''(x) + 2\lambda
\bigl(s(\infty )-s(-y_0) \bigr)^2 h_{\lambda}(x) g
\bigl(x \bigl(s(\infty)-s(-y_0) \bigr)+s(-y_0)
\bigr)=0,&\nonumber\\[-8pt]\\[-8pt]
 & h_{\lambda
}'(-\infty)=0,\qquad  h_{\lambda}(0)=1,&\nonumber
\end{eqnarray}
but in fact, setting, for any $t \in\R$,
%
\begin{equation}
f_{\lambda}(t) = {h_{\lambda} \biggl( \frac
{s(-t)-s(-y_0)}{s(\infty)-s(-y_0)} \biggr)}\Bigl/{h_{\lambda} \biggl(\frac
{-s(-y_0)}{s(\infty)-s(-y_0)}  \biggr)}
\end{equation}
allows to exactly recover the result of the corollary.
\end{remark}

\subsection{Proof of Theorem \texorpdfstring{\protect\ref{th:lim3}}{2.3}}
\label{sec:toinfty}
The framework is now the one of Theorem~\ref{th:lim3}, that is
Assumptions \ref{assume:V}, \ref{assume:start2} are in force, in particular
for odd $d$ we have
$y_0^{(\varepsilon)}\to-\infty$, which represents a novelty with
respect to what
we have done
up to now. Therefore we start with a result that addresses this issue.

\begin{lemma}
\label{lem:convy0}
The limits in distribution of $T_{y}$, as $y \to-\infty$, and of
$\tau_{a \gep^{-2/d}}(Y)$, as $\gep\searrow0$, exist and coincide.
\end{lemma}

\begin{pf}
Write $Y^y$ for the process $Y$ satifying (\ref{eq:rescdiff}) and
$Y_0=y$, and
simply observe that for any $y$ such that $y_0^{(\varepsilon)}\le y\le
a\gep^{-2/d}$
%
\begin{equation}
\label{eq:decomptau} \tau_{a \gep^{-2/d}} \bigl(Y^{y_0^{(\varepsilon)}} \bigr) = \tau
_{y} \bigl(Y^{y_0^{(\varepsilon)}} \bigr) + \tau_{a\gep^{-2/d}}
\bigl(Y^y \bigr) .
\end{equation}
By Proposition~\ref{prop:limy_0}, the second term in the above sum
converges in distribution to $T_{y}$ as $\gep\searrow0$.
The proof will therefore be completed if we show that $\tau
_y(Y^{y_0^{(\varepsilon)}})$, the time taken by our rescaled diffusion
to go from
$y_0^{(\varepsilon)}$ to some fixed $y$, becomes negligible when $y
\to-\infty$,
uniformly over small $\varepsilon$.

The strategy of proof is quite straightforward: for $y$ large (and
negative), the diffusion process $Y_t$, while it stays in
$[y_0^{(\varepsilon)}, y]$,
is dominated (with overwhelming probability) by the drift term. But it
is straightforward to see that, under Assumption~\ref{assume:V},
such a non-Lipschitz drift term drives the solution from arbitrarily
far to the left to $y$ in a finite time, which can even be made
arbitrarily small by choosing $\vert y \vert$ sufficiently large. To
detail these steps let us introduce for $M>0$ the event
%
\begin{equation}
\gO_M := \bigl\{ \bigl\vert B(t)\bigr\vert\le M \mbox{ for every } t\in
[0,1] \bigr\} ,
\end{equation}
where $B(\cdot)$ is the Brownian motion that drives $Y$. Of course the
probability of $\gO_M$ tends to one as $M$ becomes large.
Let us work with $M=\vert y \vert/4$ and let us assume that $\gO_M$
is verified and that $\vert y \vert\ge A$, $A$
given in Assumption~\ref{assume:V}.
As a first step let us introduce $U_t$ strong solution of $\mathrm{d} U_t=
\psi(\vert U_t\vert) \dd t + \dd B_t$
for $U_0=y_0^{(\varepsilon)}$: note that, in view of Assumption~\ref
{assume:V},
there is no loss of generality
in assuming that $\psi(\cdot)$ is not simply continuous, but smooth.
Moreover in view of Lemma~\ref{th:local} the shape of $V\me(y)$ for
$y\le b \gep^{-2/d}$ is inessential, so we may
as well assume that $V'\me(y) \le-\psi(\vert y\vert)$ also for
$y\le b \gep^{-2/d}$.
Then one directly checks by analyzing
$\mathrm{d}(Y_t-U_t)$ that $Y_t \ge U_t$ for every {$t\le\zeta:=\inf
\{t\dvt  U_t \ge y\}$}, because $y\le-A$. Therefore, $\zeta$ gives
an upper bound on $\tau_{y}(Y^{y_0^{(\varepsilon)}})$. Let us then
focus on
$U_\cdot$ and let us set
$Z_t:= U_t-B_t$, so $Z_\cdot$ is differentiable and

%
\begin{equation}
\label{eq:akd1} \frac{\mathrm{d}}{\mathrm{d}t} Z_t = \psi \bigl( \vert
Z_t+B_t\vert \bigr) \ge\psi ( -Z_t-M ) ,
\end{equation}
where the inequality holds for $t\in[0, \zeta\wedge1]$ and on $\gO
_M$: note in fact that for
$t\le\zeta$ we have $Z_t+B_t=U_t\le y <0$ so, since we have also
$t\le1$, $ \vert Z_t+B_t\vert=-Z_t-B_t\ge-Z_t-M$.
By integrating the differential inequality \eqref{eq:akd1}, we
see that on $\gO_M$
%
\begin{eqnarray}
\zeta\wedge1 &\le&\int_{y_0^{(\varepsilon)}}^{Z_{\zeta\wedge1}} \frac{\mathrm{d}
u}{\psi(-u -M)}
 \le \int_{y_0^{(\varepsilon)}}^{y+M} \frac{\mathrm{d}
u}{\psi(-u
-M)}
\nonumber
\\[-8pt]\\[-8pt]
& = & \int_{\vert y\vert/2}^{-y_0^{(\varepsilon)}-\vert y\vert/4} \frac
{\mathrm{d} u}{\psi(u)} \le\int
_{\vert y\vert/2}^{\infty} \frac{\mathrm{d} u}{\psi(u)} ,\nonumber
\end{eqnarray}
where we used an obvious change of variable and $M=|y|/4$ to get the
equality of the second line.
Since the rightmost term above can be made arbitrarily small
by choosing $\vert y \vert$ large, $\zeta$ is at most of this size
and the proof is complete.
\end{pf}

\begin{pf*}{Proof of Theorem~\ref{th:lim3}}
For even $d$, Theorem~\ref{th:lim3} follows directly from Proposition~\ref{prop:limy_0} once we take into account Remark~\ref{rk:y_0ep}.
Note that in this case $T_0$ of Proposition~\ref{prop:limy_0} is the
random variable $T_{d,\mu}$ of Theorem~\ref{th:lim3}.
For odd $d$ instead Proposition~\ref{prop:limy_0} has to be combined
with Lemma~\ref{lem:convy0}:
the random variable identified by the limit(s) in Lemma~\ref
{lem:convy0} is $T_{d,\mu}$.
\end{pf*}

We are also ready to generalize Corollary~\ref{cor:edo} to cover the
case $y_0=-\infty$ when $d$ is odd.

\begin{cor}
\label{cor:edo2}
Under the same assumptions as in Corollary~\ref{cor:edo} (in
particular recall that $\lambda\le0$), when $d$ is odd
the limits $\lim_{y \to\pm\infty}f_\gl(y)=: f_\gl(\pm\infty)
\in(0, \infty)$ exist and
%
\begin{equation}
\label{eq:coredo2} \ee \bigl[\exp (\lambda T_{d,\mu} ) \bigr] =
\frac{f_{\lambda}(-\infty)}{f_{\lambda}(\infty)} .
\end{equation}
\end{cor}

\begin{pf} The fact that $\lim_{y \to\infty}f_\gl(y)= f_\gl
(\infty) \in(0, \infty)$ has been already proven in
Corollary~\ref{cor:edo} and it is exactly from \eqref{eq:coredo} that
we restart.
Since $\gl\le0$, this is just a matter of taking $y_0\to-\infty$ in
\eqref{eq:coredo},
by using the fact that the limit of the left-hand side exists by
Proposition~\ref{prop:limy_0}.
\end{pf}

Let us conclude the section with two remarks: the first is particularly
important.

\begin{remark}
\label{rem:analytic}
For sake of conciseness, we have chosen to focus on convergence in law
of $\{\gep^{2/d} \tau_a (X)\}_{\gep>0}$ and on
properties of the limit $T_{d,\mu}$.

So far we have proven the convergence in law. In Section~\ref{sec:Schr} below, we will check that assumptions of Corollary~\ref
{cor:edo} are satisfied (see Remark~\ref{rem:dersub} below), which
ensures that
%
\begin{equation}
\label{eq:convexp} \lim_{\gep\searrow0} \bbE \bigl[\exp \bigl(\gl
\gep^{2/d} \tau_a (X) \bigr) \bigr] = \frac{f_\gl(l_d)}{f_\gl(\infty)} ,
\end{equation}
with $l_d=-\infty$ for $d$ odd and $l_d=0$ for $d$ even, but \emph{a
priori} only for $\gl\le0$. We then resort to analytic continuation
arguments to show that
%
\begin{equation}
\label{eq:coredo2.1} \Phi_{T_{d, \mu}} ( \gl ) = \ee \bigl[\exp (\lambda
T_{d,\mu} ) \bigr] = \frac{f_{\lambda}(l_d)}{f_{\lambda}(\infty)} ,
\end{equation}
holds also for $\gl>0$: precisely, it holds for every $\gl< \gl
_0=\gl_0(T_{d,\mu})$ (which was defined in the beginning of Section~\ref
{sec:mainres}). Theorem~\ref{th:lim3} states in addition that it
coincides with $\widetilde{\gl}_{0,d}$, the detailed proof of this
fact is given in Section~\ref{sec:ode} for $d$ odd, and in Section~\ref{sec:ode2} for $d$ even.

More precisely we will show that
the functions of $\gl$ appearing in the denominator and numerator
of the rightmost side of \eqref{eq:coredo2.1} can both be extended to
the whole complex plane as entire functions (e.g., in the odd
case, see (\ref{eq:laplace--g}), (\ref{eq:Stransf2}), (\ref
{eq:subentire}), and Theorem~\ref{th:WKBx} below),
therefore the only singularities of the rightmost side of \eqref
{eq:coredo2.1} are poles and so it is analytic (in particular) in
$\{\gl\in\bbC\dvt  \Re(\gl) < \tilde\gl\}$; in principle $\tilde\gl
= \min\{ \Re(\gl) \dvt  f_{\lambda}(\infty)=0\}$,
but it turns out all the poles are real so $\tilde\gl$ is itself a
pole.

Of course $\tilde\gl\ge0$, but $\tilde\gl=0$ is excluded too since
$\vert\Phi_{T_{d, \mu}}(\gl)\vert\le1$ for every $\gl\in \mathrm{i} \bbR$,
and this is incompatible with the fact that the singularity of the
rightmost side of \eqref{eq:coredo2.1} is a pole.

On the other hand, by definition of $\gl_0$ one directly sees that
$\Phi_{T_{d, \mu}}(\cdot)$
is analytic in $\{\gl\in\bbC\dvt  \Re(\gl) < \gl_0\}$.
So
\eqref{eq:coredo2.1} holds for $\gl< \tilde\gl\wedge\gl_0$ and we
are left with proving that
$\tilde\gl= \gl_0$.
For this, we observe that $\tilde\gl< \gl_0$ means that
$\Phi_{T_{d, \mu}}(\gl)$ blows up as $\gl\nearrow\tilde\gl$,
which is impossible by the definition
of $\gl_0$. On the other hand, $\tilde\gl> \gl_0$ cannot hold
either: to show this assume $\tilde\gl> \gl_0$.
If $\Phi_{T_{d, \mu}}(\gl_0)=\infty$, the contradiction is
immediate. Let us therefore assume also that $\Phi_{T_{d, \mu}}(\gl
_0)<\infty$ and
observe that, since
the expression in
the rightmost term in \eqref{eq:coredo2.1}
has radius of convergence $\tilde\gl-\gl_0>0$ at $\gl_0$, the
$n$th derivative in $\gl$ of this term
is $\mathrm{O}(\gd^{-n})$, for any $\gd<(\tilde\gl-\gl_0)$. But, with the notation
$T=T_{d, \mu}$, this implies
$\bbE[ T^n \exp( \gl_0 T)] =\mathrm{O}(\gd^{-n} n!)$, from which one
directly extracts that $\bbE[ \exp( \gl T)] < \infty$
for $\gl< \gl_0+ \gd$, which contradicts the assumption.
Therefore $\tilde\gl=\gl_0$ and, for any $\gl\in\bbC$ such that
$\Re(\gl) $ is smaller than this value, \eqref{eq:coredo2.1} holds.
\end{remark}

\begin{remark}
\label{rem:exp-int}
It is natural to wonder about the validity of \eqref{eq:convexp} for
some $\gl>0$,
possibly for all $\gl< \gl_0(T_{d, \mu})$, as
Theorem~\ref{th:Phiright} may suggest. This however, it is in general
false: in fact, it may be that
$\bbE [\exp (\gl\gep^{2/d} \tau_a (X) )
 ] = \infty$ for every $\gl>0$ and every $\gep>0$, see Figure~\ref{fig:esc-left}.
%
\begin{figure}[b]

\includegraphics{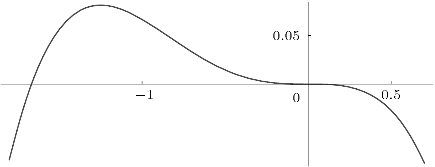}

\caption{The choice $U_{0, \gep}(x)= -\frac{1}6 x^3 - \frac{1}{10}x^4$,
independent of $\gep$, fulfills
Assumption~\protect\ref{assume:V} with rescaled potential $V_0(y)=-\frac{1}6
y^3$. Our results apply to the diffusion $X_\cdot$, cf. \protect\eqref{eq:diffusion},
with $X_0=-1$, and we know
that
$\gep^{2/3}\tau_{a, \gep}(X)$ (cf. \protect\eqref{eq:tau}), any $a>0$,
converges as $\gep\searrow0$ to $T_{3,0}$
and this limit random variable has exponential moments (Theorem~\protect\ref
{th:Phileft}). However, since $X_\cdot$ reaches $-\infty$
with positive probability, hence
$\bbP(\tau_{a, \gep}(X)=\infty)>0$ for every $\gep>0$ and therefore
$\bbE [\exp (\gl\gep^{2/d} \tau_{a, \gep} (X) )
 ] = \infty$ for every $\gep>0, \gl>0$.}\label{fig:esc-left}
\end{figure}

Other, more subtle phenomena may happen and they are notably connected
to the fact that we are making assumptions
only on the limit behavior of $V_{\mu, \gep}$ and not on $V'_{\mu,
\gep}$. We make the choice not to go
further toward this direction and we just stress that stronger
assumptions on $U_{\mu, \gep}$, or on $V_{\mu, \gep}$,
are needed to establish \eqref{eq:convexp}, that is the convergence of
the moment generating function, for $\gl\in(0, \gl_0)$.
\end{remark}
%


\section{ODE analysis: The saddle node case}
\label{sec:ode}

For sake of conciseness, as we explained in Section~\ref{sec:main},
we restrict to $d=3$, even if what we present can be generalized to arbitrary
odd $d$ (see Remark~\ref{rem:replace} below -- only for Corollary~\ref{cor:sharp} below does the generalization become cumbersome).
Corollary~\ref{cor:edo} and Corollary~\ref{cor:edo2} tell us that we
are interested in a specific solution to
\eqref{eq:edofgl}: in this section we show the existence of such a
solution and prove a number of quantitative results,
that yield a proof of Theorem~\ref{th:Phiright} and Theorem~\ref{th:Phileft}
for the saddle node case.

\subsection{Schr\"odinger equation: Mapping, basic facts}
\label{sec:Schr}
For conformity with the ODE literature \cite{cf:CL,cf:Sibuya,cf:T1},
we set $\eta= 2\lambda$ and focus, for $\mu\in\R$ and $\eta\in
\bbC$, on solutions $g=g_\eta$ to
%
\begin{equation}
\label{eq:g} g_\eta''(x) +
\bigl(x^2- \mu \bigr) g_\eta'(x) +\eta
g_\eta(x) = 0 ,
\end{equation}
which of course is just a rewriting of \eqref{eq:edofgl}.
More precisely we look for a solution $g_{\eta}$ to (\ref{eq:g}) such
that $g'_{\eta}(-\infty)=0$. We will actually look for a solution
$g_{\eta}$
that is bounded at $-\infty$, and we will then see that such a
solution satisfies $g'_{\eta}(-\infty)=0$. In any case, this boundedness
requirement can (and, we will see, does) determine $g_\eta(\cdot)$,
but only up to a multiplicative constant.
This is a side issue
for the moment since (for $\eta\le0$, recall Corollary~\ref{cor:edo2})
%
\begin{equation}
\label{eq:laplace--g} \Phi_{T_{3, \mu}}(\lambda) = \Phi_{T_{3, \mu
}}(\eta/2) =
\ee \biggl[\exp \biggl(\frac{\eta}{2} T_{3, \mu} \biggr) \biggr] =
\frac{g_{\eta}(-\infty)}{g_{\eta}(\infty)} ,
\end{equation}
whenever such a solution $g_{\eta}$ exists.

Note that if we set
%
\begin{equation}
\label{eq:Stransf} g_\eta(x) := \exp \bigl(V_\mu(x)
\bigr)u_\eta(x) = \exp \biggl(\mu\frac{x}2 -
\frac{x^3}6 \biggr)u_\eta(x) ,
\end{equation}
then $u=u_\eta$ -- when the context is clear we drop the subscript
$\eta$ from the notation -- solves
the Schr\"odinger equation
%
\begin{equation}
\label{eq:S} u''(x)- Q_{\eta}(x) u(x) =
0
\end{equation}
with
%
\begin{equation}
\label{eq:Qeta9} Q_{\eta}(x) := \bigl(V_{\mu}'
\bigr)^2(x)-V_{\mu}''(x) - \eta=
x + \tfrac{1}{4} \bigl(x^2-\mu \bigr)^2 - \eta =:
q(x) - \eta.
\end{equation}
The ODE (\ref{eq:S}) has been intensively studied, see in particular
\cite{cf:CL,cf:T1}
and the more recent \cite{cf:Sibuya}. Therefore we start by reminding
some of the results which we need in our study, and then we state our results.
%
\begin{figure}

\includegraphics{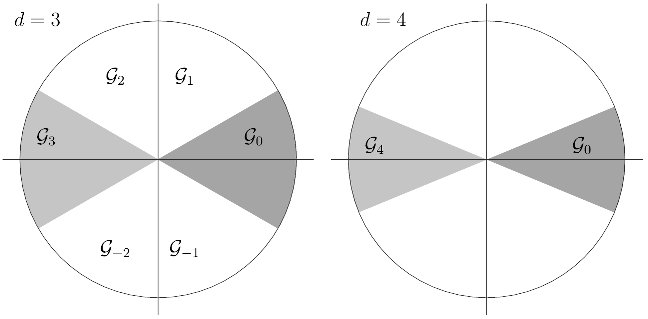}

\caption{The sectors of $\bbC$ relevant for the study of \protect\eqref
{eq:S}, for the cases $d=3$ and $d=4$. It is not difficult
to understand why these sectors are singled out: for $\vert x \vert$
large $Q_\eta(x)\sim c_d x^m$ ($c_d>0$ and
$m=2d-2$). Via a quick WKB analysis (made for example, at an informal
level like in \cite{cf:BO}, pages 486--487) one readily realizes
that the leading contribution to the asymptotic behavior for $\vert
x\vert$ large is $\exp(c'_d x^{(m+2)/2})$
(with $c'_d=2c_d/(m+2)$), at least if $\Re(x^{(m+2)/2})\neq0$. If we
consider $x=r\exp(\mathrm{i} \theta)$, $r>0$, $\Re(x^{(m+2)/2})=0$
if $\theta=\uppi(1+2k)/(2d)$, for some integer $k$. So the asymptotic
limits switch behavior when \textsl{crossing}
the rays corresponding to the \textsl{critical} angles separating the
sectors.}\label{fig:sectors}
\end{figure}

For $k \in\{-2,-1,0,1,2,3\}$ define the open sector $\mathcal{G}_k =
\{ x \in\mathbb{C} \dvt  \llvert \arg(x) - \frac{k \uppi}{3}\rrvert  <
\frac{\uppi}{6} \}$.z
Two of these sectors are of particular interest: $\mathcal{G}_0$
contains $(0, \infty)$, while $\cG_3$ contains $(-\infty,0)$, see
Figure~\ref{fig:sectors} and its caption for further explanation.
By \cite{cf:Sibuya}, Th. 6.1 and Th. 7.1, for each $\cG_k$, there
exists a unique solution to (\ref{eq:S}), called the \emph
{subdominant solution in the sector $\cG_k$}, and denoted $u_{\mathrm
{sub},\cG_k,\eta}=u_{\mathrm{sub},\cG_k}$, which, up to a
multiplicative constant, is characterized by
%
\begin{equation}
\lim_{\stackrel{|x| \to\infty, x \in\cG_k}{\arg(x) \mathrm{\ fixed}}} u_{\mathrm{sub},\cG_k}(x) = 0.
\end{equation}
The solution $u_{\mathrm{sub},\cG_3}$ which is subdominant along the
negative real line, will be of particular interest for us and it is
fully (i.e., not just up to a multiplicative constant)
characterized by
%
\begin{equation}
\label{eq:G3sub} \lim_{\stackrel{|x| \to\infty, x \in\cG_3}{\arg(x) \mathrm{\ fixed}}} \exp \biggl( -\frac{x^3}{6} +
\frac{\mu}{2} x \biggr) u_{\mathrm{sub},\cG_3}(x) = 1 .
\end{equation}
So we choose to work with the solution $g(\cdot)=g_\eta(\cdot)$
%
\begin{equation}
\label{eq:Stransf2} g_\eta(x) = \exp \biggl( - \frac{x^3}6 + \mu
\frac{x}2 \biggr)u_{\mathrm{sub},\cG_3}(x)\stackrel{x \to-\infty}\sim1 ,
\end{equation}
where the asymptotic equivalence is of course a particular case of
\eqref{eq:G3sub}.
So, with this choice, namely $g_\eta(-\infty)=1$, the expression
\eqref{eq:laplace--g}
becomes somewhat simpler and, above all, we can directly apply results in
\cite{cf:Sibuya} that grant analyticity. In fact
Theorem~6.1 and Theorem~7.1 in \cite{cf:Sibuya} ensure that for any $k
\in\{-2,-1,0,1,2,3\}$,
%
\begin{equation}
\label{eq:subentire} u_{\mathrm{sub},\cG_k}(x) \mbox{ is entire with respect to } (\mu ,
\eta,x) .
\end{equation}

\begin{remark}
\label{rem:dersub}
In \cite{cf:Sibuya}, Th. 6.1, and explanations following it one finds
an expansion to all orders of $u'_{\mathrm{sub},\cG_3}$ in the same
limit as in
\eqref{eq:G3sub}. More precisely, as $|x| \to\infty, x \in\cG_3,
\arg(x) \mathrm{\ fixed}$,
%
\begin{equation}
\label{eq:dersubG3} u'_{\mathrm{sub},\cG_3}(x) = \biggl(
\frac{x^2-\mu}{2}+\mathrm{o}(1) \biggr) \exp \biggl( \frac{x^3}{6} -
\frac{\mu}{2} x \biggr) .
\end{equation}
Observe that one formally obtains \eqref{eq:dersubG3} by
differentiating \eqref{eq:G3sub}.
This actually says that $
f_\gl(\cdot)=g_{2\gl}(\cdot)$ fulfills the assumptions in Corollary~\ref{cor:edo}. This settles the issue:
$g_\eta$ is the solution of \eqref{eq:g} we are interested in.
\end{remark}

Note that the uniqueness (this time, up to a multiplicative constant)
of $u_{\mathrm{sub},\cG_k}$ implies that any solution of (\ref
{eq:S}) which is not a multiple of $u_{\mathrm{sub},\cG_k}$
does not go to $0$ when $|x| \to\infty, x \in\cG_k, \arg(x)$
fixed. But more than that is true (see \cite{cf:Sibuya}, page 19):
any solution of (\ref{eq:S}) which is not a multiple of $u_{\mathrm
{sub},\cG_k}$ is \textsl{dominant} in $\cG_k$, that is it tends to
$\infty$
when $|x| \to\infty, x \in\cG_k, \arg(x)$ fixed. We can even be
more precise: still from \cite{cf:Sibuya}, page 19, we learn that
if a solution is subdominant in $\cG_k$, it is dominant in the
neighboring sectors $\cG_{k\pm1}$.
Therefore, by their distinct behavior at infinity, $u_{\mathrm
{sub},\cG_k}$ and (say) $u_{\mathrm{sub},\cG_{k+1}}$ are independent
and form a basis (of course this implies that any solution that is not
subdominant in a sector is necessarily dominant).
%

There is a priori no canonical choice of a dominant solution in a
sector, but it is useful to make one and ours is
$u_{\mathrm{dom},\cG_k}:=u_{\mathrm{sub},\cG_{k+1}}$. This has the
advantage that $u_{\mathrm{dom},\cG_k, \eta}(x)
=u_{\mathrm{dom},\cG_k }(x)$
is entire in $\eta$ and $x$. Of course we have
%
\begin{equation}
\label{eq:coord} u_{\mathrm{sub},\cG_3} = a u_{\mathrm{sub},\cG_0} + b u_{\mathrm
{dom},\cG_0} ,
\end{equation}
for a unique choice of $a$ and $b$, once $\eta$ is fixed,
and with such a choice we have also that $a$ and $b$ are entire
functions of $\eta$.
This follows simply from the fact that $a$ and $b$ are determined by
%
\begin{equation}
\cases{ u_{\mathrm{sub},\cG_3}(0) = a u_{\mathrm{sub},\cG_0}(0) + b
u_{\mathrm{dom},\cG_0}(0) ,
\cr
u'_{\mathrm{sub},\cG_3}(0) = a
u'_{\mathrm{sub},\cG_0}(0) + b u'_{\mathrm{dom},\cG_0}(0) , }
\end{equation}
which of course is solvable since the Wronskian matrix of $\{u_{\mathrm
{sub},\cG_0}(0),u_{\mathrm{dom},\cG_0}(0)\}$
is invertible.
In Section~\ref{sec:proofWKB}, we will provide the
precise asymptotic behaviour of $u_{\mathrm{dom},\cG_0}$ at $+\infty
$ and this will be a key step
to establish the analyticity of $\eta\mapsto g_\eta(\infty)$.

\subsection{The spectrum of the Schr\"odinger operator} \label{sec:spectrum}
A natural question to ask is wether one can indeed have $b = 0$, that
is, if it can
be that
%
\begin{equation}
\label{eq:uL2} u_{\mathrm{sub},\cG_3} \propto u_{\mathrm{sub},\cG_0} .
\end{equation}
For this note, the sharp asymptotic behavior of all the subdominant
solutions is known. In fact
\cite{cf:Sibuya}, Th. 6.1, Th. 7.1,
in analogy with \eqref{eq:G3sub}
%
\begin{equation}
\label{eq:G0sub} \lim_{\stackrel{|x| \to\infty, x \in\cG_0}{\arg(x) \mathrm{
fixed}}} x^2 \exp \biggl(
\frac{x^3}{6} - \frac{\mu}{2} x \biggr) u_{\mathrm{sub},\cG_0}(x) = 1 ,
\end{equation}
in which of course we have made a precise choice of the multiplicative constant,
but the key point is that \eqref{eq:G0sub}, coupled of course with
\eqref{eq:G3sub}, implies that if \eqref{eq:uL2} is satisfied, then
we have a solution in $\mathbb{L}^2(\bbR)$. So we have found an
eigenfunction of the differential operator $u \mapsto-u'' +q u$
with eigenvalue $\eta$
(recall \eqref{eq:S} and note the operator $L$
\eqref{eq:Schr} is just $1/2$ times the operator we are considering
here). On the other hand, since we have seen
that $\{u_{\mathrm{sub},\cG_0, \eta}, u_{\mathrm{dom},\cG_0, \eta
}\}$ is a basis, $\eta$ is an eigenvalue in $\mathbb{L}^2(\bbR)$
if and only if \eqref{eq:uL2} holds.

Therefore the question we have just raised is about the spectrum of the
Schr\"odinger operator.
It is well known -- see, for example, \cite{cf:CL}, Ch. 8, or \cite{cf:Sibuya}, Ch. 7 -- that the spectrum is constituted by an infinite
sequence of eigenvalues $\eta_0< \eta_1 < \cdots\,$,
$\spect:=\{\eta_0, \eta_1, \ldots\}$: the eigenfunction $u_k$ of
$\eta_k$ is $u_{\mathrm{sub},\cG_3, \eta}$, or $u_{\mathrm
{sub},\cG_0, \eta}$.
Moreover, as it is explained in detail in \cite{cf:CL}, Ch. 8, and
\cite{cf:Sibuya}, Ch. 7, for $\eta$ on the real axis
and fixed boundary conditions the number of the zeros is decreasing in
$\eta$ and the location of the zeros is a continuous function of $\eta
$. In particular
$u_{\mathrm{sub},\cG_3, \eta}(x)$, and therefore $g_\eta(x)$, has
no zero for $\eta<\eta_0$.
Note that in our case $\eta_0>0$: in fact, \eqref{eq:Stransf}
directly implies that $g_\eta(\infty)=0$ for $\eta$ in the spectrum
and we have chosen $g_\eta(-\infty)=1$, so the right-hand side of
\eqref{eq:laplace--g} is $\infty$,
but the left-hand side is bounded by $1$ for $\eta\le0$, so $\spect
\subset(0, \infty)$.

\begin{remark}
\label{rem:1at0}
This tells us in particular that for $\eta< \eta_0$, we can certainly
redefine $g_\eta(\cdot)$
so that $g_\eta(0)=1$. With such a choice of course $g_\eta(\cdot)$
would still be analytic,
but it is not defined for all $\eta\in\bbC$, since for (infinitely
many) real values of
$\eta< \eta_0$ we have $g_\eta(0)=0$.
\end{remark}


\subsection{Analysis of the ODE: The results}
\label{sec:ODEres}
We now state the precise estimates for $u_{\mathrm{sub},\cG_3, \eta
}$, in both asymptotic limits
\begin{enumerate}[(2)]
\item[(1)] $|\eta| \to\infty$ (except along an arbitrarily small sector
containing the positive real axis), uniformly in $x$;
\item[(2)] $\eta\in\mathbb{C} \setminus\spect$ is fixed and $x \to
\infty$.
\end{enumerate}
Before stating both results, we need to precise that for any $x \in\R
$, when $Q_{\eta}(x)\in\bbC\setminus(-\infty, 0]$
(which is always the case when $\eta$ is large, and in a sector of the
complex plane that does not contain $(0, \infty)$), we define $
Q_{\eta}^{1/2}(x)$ (resp., $Q_{\eta}^{1/4}(x)$) as the only square
root of $Q_{\eta}(x)$ (resp., $Q_{\eta}^{1/2}(x)$) which satisfies
$\arg(Q_{\eta}^{1/2}(x)) \in(-\uppi/2, \uppi/2]$ (resp.,
$\arg(Q_{\eta}^{1/4}(x)) \in(-\uppi/4, \uppi/4]$). This corresponds of
course to choosing what is normally called the \textsl{principal
branch} of the square root.

\begin{theorem}
\label{th:WKBeta}
Fix $\mu\in\bbR$ and $\alpha\in(0,\uppi)$. For any $\eta\in
\mathbb{C} \setminus(0,\infty)$ there exists a solution $u$ of
\eqref{eq:S} satisfying, uniformly in $\eta$ such that $\vert
\arg(-\eta) \vert\le\alpha$,
%
\begin{equation}
\label{eq:WKBeta} \sup_{x \in\bbR} \biggl\llvert u (x)
Q_\eta^{1/4}(x) \exp \biggl(-\int_0^x
Q_\eta^{1/2} (y) \dd y \biggr) -1 \biggr\rrvert = \mathrm{O}
\bigl( | \eta|^{-3/4} \bigr) .
\end{equation}
\end{theorem}

Observe that the above solution is clearly subdominant for $x \to
-\infty$ and it is therefore proportional to $u_{\mathrm{sub},\cG_3,
\eta}(x)$.

We defer the proof of Theorem~\ref{th:WKBeta} to Section~\ref
{sec:proofWKB}:
from the proof it is not difficult to see that the statement actually
remains true in the case of general odd $d$,
actually the exponent in the error term can be replaced by $d/(d+1)$
(see also Remark~\ref{rem:replace}). More than that, the statement holds
also for even $d$ (Theorem~\ref{th:WKBetaquartic}).

We have similar, albeit more implicit, results when $\eta$ is fixed
and $x \to\infty$.

\begin{theorem} \label{th:WKBx}
For every $\eta\in\mathbb{C} $,
the following limit exists
%
\begin{equation}
\label{eq:WKBx} \lim_{x \to+\infty} u_{\mathrm{sub},\cG_3}(x) \exp \biggl(-
\frac{x^3}{6} + \mu\frac
{x}{2} \biggr) =: \ell(\eta) \in\bbC.
\end{equation}
Moreover $\ell(\cdot)$ is entire, $\ell(\eta)=0$ if and only if
$\eta\in\spect$ and we have $\ell(\eta)>0$ for
$\eta< \eta_0$.
\end{theorem}

We defer the proof of Theorem~\ref{th:WKBx} to Section~\ref{sec:proofWKB}.
There are several consequences to the theorems we just stated,
corresponding to the following three corollaries. The proofs of these
three corollaries are deferred to the end of this section, Section~\ref{sec:proofcorollaries}.

\begin{cor}
\label{cor:lambda0}
Recall \eqref{eq:lambda0} and \eqref{eq:lim} for the definition of
$\lambda_0=
\lambda_0(T_{3,\mu})$ and that $\eta_0$ is the smallest
eigenvalue for \eqref{eq:S}.
Then
%
\begin{equation}
\label{eq:eta0lambda0} \lambda_0 = \frac{\eta_0}{2} .
\end{equation}
Moreover, $\Phi_{T_{3,\mu}}(\gl)$ has a simple pole at $\gl_0$ and
the value of the residue $-C_{3, \mu}$ (cf. Theorem~\ref
{th:Phiright}) is
%
\begin{equation}
\label{eq:residue3} -\frac{c_{3, \mu}}{2\int_\bbR u^2_{\mathrm{sub},\cG_3, \eta_0}
(x)\dd x} ,
\end{equation}
where $c_{3, \mu}= \lim_{x \to\infty} x^2 \exp(x^3/6-\mu x/2)
u_{\mathrm{sub},\cG_3, \eta_0}(x)\in(0, \infty)$.
\end{cor}


We are further able to deduce, using Theorem~\ref{th:WKBeta}, the
asymptotic behavior of the ratio $g_\eta(-\infty)/g_\eta(\infty)$,
which with our normalization choice reduces to $1/g_\eta(\infty)$, as
$\vert\eta\vert\to\infty$ along any ray that is not the positive
real axis.

\begin{cor}
\label{cor:firstorder}
Fix $\beta_0 \in(0,\uppi)$, uniformly in $\gb\in[-\gb_0, \gb_0]$
we have that
%
\begin{equation}
\label{eq:firstorder0} \limtwo{|\eta| \to\infty} { \arg(-\eta) \to\beta}
\frac{1}{\vert
\eta\vert^{3/4}}\log \bigl(\Phi_{T_{3, \mu}}(\eta/2) \bigr) = - \limtwo {|
\eta| \to\infty} { \arg(-\eta) \to\beta} \frac{1}{\vert\eta
\vert^{3/4}}\log g_\eta(
\infty) = - c(\gb) ,
\end{equation}
%
where
%
\begin{eqnarray}
\label{eq:a_0} c(\beta) &:=& \int_{-\infty}^{\infty}
f_0^{(\beta)}(y) \dd y, \quad \mbox {and} \nonumber\\[-8pt]\\[-8pt]
 f_0^{(\beta)}(y)
&= &\frac{1}{\sqrt{2}} \sqrt{\frac
{y^4}{4}+\cos(\beta) + \sqrt{1+\cos(\beta)
\frac{y^4}{2} + \frac
{y^8}{16}}} - \frac{y^2}{2}.\nonumber
\end{eqnarray}
\end{cor}

The fact that $\vert c(\gb)\vert< \infty$ follows from
$f_0^{(\gb)}(y) =\cos(\gb)/y^2 +\mathrm{O}(1/ y^6)$. Note that, for every $y$,
$f_0^{(\beta)}(y)$ increases when $\cos(\gb)$ increases, so the even function
$\gb\mapsto c(\gb)$ decreases for $\gb\in[0, \uppi)$.

By specializing to the negative real axis $\eta<0$ (so $\beta=0$), we
are able to give asymptotic
results for $g_\eta(\infty)$, that are much sharper than the ones in
Corollary~\ref{cor:firstorder}:

\begin{cor}
\label{cor:sharp}
For $\eta\to- \infty$
%
\begin{equation}
\label{eq:etareal-sharp} \Phi_{T_{3, \mu}}(\eta/2) = \frac{1}{g_\eta(\infty)} = \bigl( 1+
\mathrm{O} \bigl(\vert\eta\vert ^{-1/4} \bigr) \bigr)\exp \bigl( - (-
\eta)^{3/4} F_0 - (-\eta )^{1/4} F_1
\bigr),
\end{equation}
where for $i=0,1$, $F_i = \int_{-\infty}^{\infty} f_i(y) \dd y$, and
%
\begin{eqnarray}
f_0(y) &=& \sqrt{1+\frac{y^4}{4}} -
\frac{y^2}{2},\nonumber
\\[-8pt]\\[-8pt]
f_1(y) &=& \frac{\mu f_0(y)}{\sqrt{y^4+4}} = \frac{1}{2}\mu \biggl(1-
\frac{y^2}{\sqrt{y^4+4}} \biggr),\nonumber %
\end{eqnarray}
that is,
$F_0= 3 \Gamma (-(3/4) )^2/ (8 \sqrt{2 \uppi} )
=\,$3.496$\ldots$ and $F_1= \mu\sqrt{2/\uppi} \Gamma(3/4)^2=\,$1.198$\ldots\,$.
\end{cor}


\subsection{Proof of Theorems \texorpdfstring{\protect\ref{th:Phiright}}{2.4}, \texorpdfstring{\protect\ref{th:Phileft}}{2.5}
and of Proposition \texorpdfstring{\protect\ref{th:regular}}{2.7} ($d$ odd)}
The proofs of the two theorems is easily disposed by referring
to some of the previous statements.
The proof of Proposition~\ref{th:regular} is instead going to require
some work.

\begin{pf*}{Proof of Theorem~\ref{th:Phiright}} By recalling Remark~\ref{rem:analytic}, one sees
that Theorem~\ref{th:Phiright} is just a restatement of Corollary~\ref{cor:lambda0}.
\end{pf*}

\begin{pf*}{Proof of
Theorem~\ref{th:Phileft}, formula
\eqref{eq:Phi3.1}}
This time what we want is just a
restatement of Corollary~\ref{cor:sharp}.
\end{pf*}

\begin{pf*}{Proof of Proposition~\ref{th:regular}}
In the whole proof of the proposition,
Corollary~\ref{cor:firstorder}, which
yields the leading asymptotic behavior of $\Phi_{T_{3, \mu}}(\gl)$
for $\vert\gl\vert\nearrow\infty$, along any ray in the complex plane,
except the positive semi-axis, plays a central role. Let us therefore
start by considering the two \textsl{standard rays}:
the negative semi-axis $\eta<0$ (Laplace transform with real argument)
and $\eta\in\Im$ (characteristic function), that is respectively,
$\gb=0$
and $\gb= \pm\uppi/2$. While the first case can be seen just as a warm
up (since it is superseded by the sharp results in
Corollary~\ref{cor:sharp}), the second case actually establishes the
validity of
\eqref{eq:forfT} for $d=3$. In these two cases, we have
%
\begin{eqnarray}
f_0^{(0)}(y) &= &\sqrt{1+\frac{y^4}{4}} -
\frac{y^2}{2} ,\nonumber
\\[-8pt]\\[-8pt]
f_0^{(\uppi/2)}(y) & =& f_0^{(-\uppi/2)}(y) =
\frac{1}{\sqrt{2}} \sqrt {\frac{y^4}{4} + \sqrt{1+ \frac{y^8}{16}}} -
\frac{y^2}{2} ,\nonumber %
\end{eqnarray}
and note that these two functions are positive: one can actually
directly check that $f_0^{(\gb)}(\cdot)>0$
if and only if $\vert\gb\vert\le\uppi/2$. Therefore $c(\gb)>0$ for
these values of $\gb\in[-\uppi/2, \uppi/2]$, and of course also
in some open interval containing $ [-\uppi/2, \uppi/2]$.
Incidentally we can compute
%
\begin{equation}
\label{eq:c0} c(0) = \frac{3 \Gamma (-(3/4) )^2}{8 \sqrt{2 \uppi}} = 3.49607\!\ldots
\end{equation}
which of course coincides with the quantity $F_0$ in Corollary~\ref{cor:sharp},
and
%
\begin{equation}
\label{eq:cpmpi} c ( \pm{\uppi}/2 ) = \frac{\sqrt{2 \uppi} \Gamma(1/4)
 (\cos(\uppi/8)- \sin(\uppi/8) )}{3 \Gamma(3/4)} = 1.33789\!\ldots.
\end{equation}

We are now ready to look at other rays.
For this note that Corollary~\ref{cor:firstorder} implies that for
every $\beta_0 \in(0,\uppi)$ and every $c<\inf_{\vert\gb\vert\le
\gb_0}c(\gb)$
there exists $C>0$ such that for every $\gl\notin\{z\dvt  \Re(z) \ge
0, \vert\mathrm{arg}(z)\vert\le\uppi- \gb_0\}$
%
\begin{equation}
\label{eq:unifboundC} \bigl\llvert \Phi_T(\gl) \bigr\rrvert \le C \exp
\bigl(-c \vert\gl\vert ^\nu \bigr) ,
\end{equation}
where $T=T_{3, \mu}$ and $\nu=3/4$. As we will see in a moment, this
estimate is relevant for us as long as $c>0$:
since $c(\cdot)$ is even and (strictly) decreasing for $\gb\in[0,
\uppi)$ we define $\gb_{\mathrm{max}}\in(0, \uppi)$ such that
%
\begin{equation}
c(\gb_{\mathrm{max}}) = 0 .
\end{equation}
Note that
\eqref{eq:cpmpi} tells us that $\gb_{\mathrm{max}}> \uppi/2$.

The following lemma outlines how the bound \eqref{eq:unifboundC}
yields existence and strong regularity result
on
the density $\dens_T(\cdot)$ of the positive random variable $T$:

\begin{lemma}
\label{th:reg}
Assume that a random variable $T$ has exponential moment generating function
$\Phi_T(\cdot)$ which can be analytically extended beyond the obvious
analyticity domain
$\{ z\dvt  \Re(z)<0\}$ to an open domain that contains the complement of
the cone
$S:=
\{z\dvt  \Re(z) >0, \vert\mathrm{arg}(z)\vert< \uppi-\gb\}$, for some
$\gb\in(\uppi/2, \uppi)$
and assume that \eqref{eq:unifboundC} holds in $S^\complement$ for
some positive constants $C$ and $c$
and for some $ \nu\in(0, 1)$. Then
$\dens_T(\cdot)$ exists and it is analytic in the cone
$\{z\dvt  \Re(z) >0, \vert\mathrm{arg}(z)\vert< \gb-\uppi/2\}$.
\end{lemma}

Lemma~\ref{th:reg} (proven below) and \eqref{eq:unifboundC} directly
yield Proposition~\ref{th:regular}, except for
\eqref{eq:f_T-sharp}. For \eqref{eq:f_T-sharp}, we recall that for a
positive random variable $X$ the expression
$\psi_X(s):=\bbE[X^{s-1}]$ is called Mellin transform of (the law of) $X$.
The fundamental domain of the Mellin transform is the open strip $\{
s\in\bbC\dvt  \ga_- <\Re(s)<\ga_+ \}$, with
$\ga_+$, respectively $\ga_-$, the largest, respectively smallest,
value such that
$\bbE[X^{s-1}]<\infty$ for every $s\in(\ga_-, \ga_+)$. Of course
$\psi_X(\cdot)$ is well defined and analytic in its fundamental domain.
Therefore if $X=\exp(T)$ (of course $T=T_{3, \mu}$), then
$\Phi_T(\gl)= \psi_{\exp(T)}(\gl+1)=: \psi(\gl+1)$ and one
easily sees than
$\ga_-=-\infty$ as well as $\ga_+= \gl_0(T)$. We now appeal to a
Tauberian Theorem in the realm of Mellin transforms,
precisely to \cite{cf:FGD}, Th. 4, part (ii) (paying attention to a
misprint in the last formula of the statement: $(\log x)^k$
has to be corrected to $(\log x)^{k-1}$). The condition to apply this
statement, that is meromorphic continuation of
$\psi(\cdot)$ on a strip larger (to the right) than the fundamental
one (in our case $\psi(\cdot)$ can be meromorphically extended to
the whole of $\bbC$) and
$\psi(s)=\mathrm{O}(\vert s\vert^r)$ for
some $r>1$ and as $s$ that tends to infinity in a suitable strip around
$\ga_+$ of the form
$\{s\in\bbC\dvt  \ga_+ - \gd<\Re(s)<\ga_+ + \gd' \}$, with $\gd$
and $\gd'$ positive numbers. In our case
such a result is (largely!) achieved by \eqref{eq:unifboundC}, and
actually with arbitrary $\gd$.
What is quantitatively relevant in applying \cite{cf:FGD}, Th. 4, part
(ii), is $\gd'$ and we choose it to be so that
$\{s\in\bbC\dvt  \ga_+ - \gd<\Re(s)<\ga_+ + \gd' \}$ contains only
the pole at $s= \gl_0(T)+1$, that, by
Theorem~\ref{th:Phiright} or, equivalently, by Corollary~\ref{cor:lambda0}, is a simple pole of which we know the residue
$C:=C_{3, \mu}$. The net result
is that for the density $\dens_X(\cdot)$ of $X$ we have for $x \to
\infty$
%
\begin{equation}
\label{eq:fromFGD} \dens_X(x) = C x^{-\gl_0(T) -1} + \mathrm{O}
\bigl( x^{-b-1} \bigr) ,
\end{equation}
with $b> \gl_0(T)$ and smaller than the second smallest eigenvalue of
the Schr\"odinger operator we are dealing with.
At this point, we just use the elementary relation
%
\begin{equation}
\dens_T(t) = \exp(t) \dens_X \bigl( \exp(t) \bigr) ,
\end{equation}
and we obtain \eqref{eq:f_T-sharp}. This completes the proof of
Proposition~\ref{th:regular}.
\end{pf*}

\begin{remark}\label{rk:expansion}
Th. 4, part (ii) of \cite{cf:FGD} gives a formula for the asymptotic
behavior of the density in terms
of the residues of all the poles in the region of meromorphic extension
of $\psi(\cdot)$. By generalizing
Corollary~\ref{cor:lambda0} to deal also with more than just the bottom of the
spectrum it is certainly possible to get to an asymptotic formula
for which each term in the expansion corresponds to a point in the spectrum.
\end{remark}

\begin{pf*}{Proof of Lemma~\ref{th:reg}}
The existence of $\dens_T(\cdot)\in C^\infty$ is a direct
consequence of the decay of the
characteristic function $\gp_T(s)= \Phi_T(\mathrm{i}s)= \mathrm{O}( \exp(-c \vert s
\vert^\nu))$, $ \vert s \vert\to\infty$:
in fact the inversion formula
$2 \uppi\dens_T(t)= \int_{-\infty}^\infty\exp(- \mathrm{i} t s) \gp_T(s)
\dd s$
%
\begin{figure}

\includegraphics{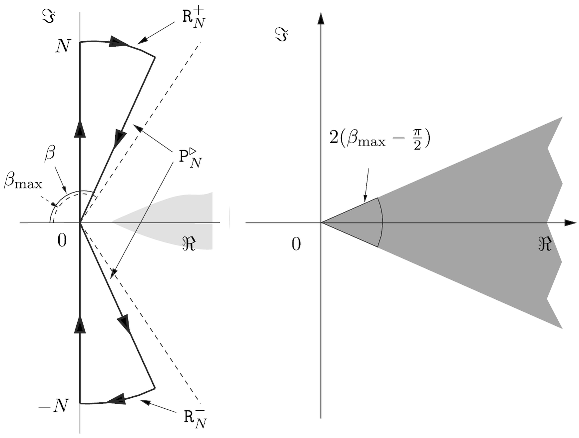}

\caption{On the left the integration path $\mathtt{C}_N$, which is the union of
the straight path going from $-\mathrm{i}N$ to $+\mathrm{i}N$,
the two arcs $\mathtt{R}^\pm_N$
and $\mathtt{P}^{\triangleright}_N$.
This path is contained in the domain of analyticity of the analytic
continuation of the moment generating function
$\Phi_T(\cdot)$, which is the complement of the lightly shadowed
region. Actually, in our
specific cases $\Phi_T(\cdot)$ can be extended to a meromorphic
functions that has (countably many) poles
on the positive real axis (that accumulate only at $+\infty$). The
angle $\gb$ that determines the path $\mathtt{C}_N$
can be chosen arbitrarily close to $\gb_{\mathrm{max}}$. The rays
corresponding to $\gb_{\mathrm{max}}$ in the figure,
that is the dashed rays, are the critical rays along which $\Phi_T(\gl
)$ does not have an exponential behavior with
the rate $\vert\gl\vert^\nu$: in the region
in which the path $\mathtt{C}_N$ lies the asymptotic behavior along
rays is exponentially vanishing. This makes negligible the contribution
of the arcs to the path integral, in the limit $N \to\infty$.
On the right instead we draw the cone in which $\dens_T(\cdot)$ is
analytic, as a consequence of the
bounds on $\Phi_T(\cdot)$ underlying the figure on the left.}\label{fig:contour}
\end{figure}
holds as soon as $\gp_T(\cdot)\in\bbL^1$ and, since $\gp_T(\cdot
)\in\bbL^p$ for every $p\ge1$, we have also that
$2 \uppi f^{(n)}_T(t)= \int_{-\infty}^\infty(-\mathrm{i}s)^n\exp(- \mathrm{i} t s) \gp
_T(s) \dd s$, where we have introduced the obvious notation for the
$n$th-derivative.
For the analyticity, we use Cauchy's Theorem: choose $t>0$, for $N$ positive
%
\begin{equation}
\label{eq:Cauchy} 0 = \int_{\mathtt{C}_N} \gl^n \exp(-t
\gl) \Phi_T(\gl) \dd\gl= \int_{-\mathrm{i}N}^{\mathrm{i}N}
\cdots+ \sum_{q=\pm} \int_{\mathtt{R}^q_N}
\cdots+\int_{\mathtt{P}^\triangleright_N} \cdots,
\end{equation}
where the contour $\mathtt{C}_N$ and paths $\mathtt{R}^\pm_N$ and
$\mathtt{P}^\triangleright_N$ in the complex plane
are given in Figure~\ref{fig:contour}.
The first term in the right-hand side converges, as $N \to\infty$, to
 $-\mathrm{i} (-1)^n 2 \uppi f^{(n)}_T(t)$,
whereas the fast decay at infinity of $\Phi_T(\cdot)$ immediately entails
that the contributions due to the integration along the arcs $\mathtt
{R}^\pm_N$ vanish in the same limits.
Therefore
%
\begin{equation}
\label{eq:newdern} f^{(n)}_T(t) = \lim_{N \to\infty}
\frac{1}{2 \uppi \mathrm{i}} \int_{\mathtt
{P}^\triangleleft_N} (-\gl)^n \exp(-t
\gl) \Phi_T(\gl) \dd\gl,
\end{equation}
where $\mathtt{P}^\triangleleft_N$ is the path $\mathtt
{P}^\triangleright_N$ with reversed orientation.
From \eqref{eq:newdern} and the hypotheses,
we have
%
\begin{equation}
\label{eq:newdernbound} \bigl\llvert f^{(n)}_T(t) \bigr\rrvert
\le \frac{C}{\uppi} \int_0^\infty
r^n \exp \bigl(-r t \sin(\gb ) \bigr) \dd r = \frac{C}\uppi
\frac{n!}{(t \sin(\gb))^{n+1}} .
\end{equation}
Therefore the radius of convergence of the Taylor series at $t>0$ is
(at least) $ t \sin(\beta)$,
which is the analyticity property claimed in the statement.
\end{pf*}

\begin{remark}
\label{rem:replace}
As we have already mentioned,
all the results we present here can be generalized
in a straightforward way to general $d \ge3$: let us quickly {discuss} here how the
ODE results we just presented generalize for $d$ odd.
\begin{itemize}
\item
The result of Theorem~\ref{th:WKBeta} remains true, in fact we get a
more precise estimate for greater $d$, as $\mathrm{O}(|\eta|^{-3/4})$ should be
replaced in the general case with $\mathrm{O}(|\eta|^{-d/(d+1)})$.
\item
The result of Theorem~\ref{th:WKBx} remains also valid when replacing
$-x^3/6 + \mu x/2$ with the more general $V_{\mu}(x)$ of \eqref{eq:assumeV}.
\item
The result of Corollary~\ref{cor:lambda0} also holds, with the obvious
modifications for the value of $C_{d,\mu}$.
\item
The result of Corollary~\ref{cor:firstorder} still holds when the
exponent of $|\eta|$ is replaced with\vspace*{-2pt} $\frac{1}{2} + \frac
{1}{2(d-1)} = \frac{d}{2(d-1)}$, and $f_0^{(\gb)}(y)$, with respect
to the formula in \eqref{eq:a_0}, is obtained by replacing the exponents
$2$, $4$ and $8$ with $(d-1)$, $2(d-1)$, and $4(d-1)$.
\end{itemize}
It is also interesting to observe that $f_0^{(\gb)}$ does not depend
on the parameter $\mu$.
In fact, an expression for $c(\gb)$ in terms of special
functions can be found for every $\gb$, and numerical evidence
suggests that $\gb_{\mathrm{max}}=(d-1)\uppi/d$.
\end{remark}

\subsection{Analysis of the ODE: Proof of Theorems \texorpdfstring{\protect\ref{th:WKBeta}}{4.3}
and \texorpdfstring{\protect\ref{th:WKBx}}{4.4}}
\label{sec:proofWKB}

\begin{pf*}{Proof of Theorem~\ref{th:WKBeta}}
We start by making easy observations explaining why $|\eta|$
large enough and $\arg(-\eta) \in[-\alpha,\alpha]$ simplifies our study.
Recall $0\le\alpha< \uppi$, and $Q_{\eta}, q$ have been introduced in
(\ref{eq:Qeta9}).
With our two assumptions on $\eta$ we either have $\Im(\eta) \ne0$
or $\eta<0$ and large enough so that $\eta< \min_{x \in\R} q(x) =:
q(x_0)$. Either way $Q_{\eta}(x) \in\bbC\setminus(-\infty,0]$ for
every $x \in\R$, so that both $Q_{\eta}^{1/2}, Q_{\eta}^{1/4}$ are
analytic.

Moreover, note that if $\arg(q(x_0)-\eta) \in[-(\uppi+\alpha
)/2,(\alpha+\uppi)/2]$ (which holds for $|\eta|$ large enough, in the
sector we consider), we have
%
\begin{equation}
\label{eq:minrealpart} \inf_{x \in\R} \Re \bigl(Q_{\eta}^{1/2}(x)
\bigr) = \Re \bigl(Q_{\eta}^{1/2}(x_0) \bigr) \ge
\cos \biggl(\frac{\alpha+\uppi
}{4} \biggr) \sqrt{|\eta|}. 
\end{equation}
For the imaginary part, we will only need the trivial bound
%
\begin{equation}
\label{eq:maximpart} \sup_{x \in\R} \bigl\vert\Im \bigl(Q_{\eta}^{1/2}(x)
\bigr) \bigr\vert= \sqrt{\bigl|q(x_0) + \eta\bigr|} \sin \biggl(
\frac{\alpha+\uppi}{4} \biggr) \le\sqrt{|\eta|}, 
\end{equation}
where we have taken $|\eta|$ large enough that the second inequality
above holds.

We now exploit the fact that WKB theory gives us a guess for the
asymptotic behavior of $u$ as $\eta$ tends to infinity away from
the positive axis and the next steps is writing an integral equation,
that is \eqref{eq:voc2}, for the ratio between $u$ and the WKB guess.
Of course, these steps are applications to our context of ideas taken
from the rigorous approach to WKB estimates (see, e.g.,
\cite{cf:T1} or the more recent \cite{cf:O} and references therein).

Let us then set $\xi(x):= \int_0^x Q^{1/2}(y) \dd y$. If we set
$n:= Q^{1/4} u$, we then have
%
\begin{equation}
\label{eq:n1} n''(x)-\frac{Q'}{2Q}
n' -Qn - \tilde R n = 0 ,
\end{equation}
with
%
\begin{equation}
\tilde R := \frac{Q''}{4Q}- \frac{5 (Q')^2}{16 Q^2} .
\end{equation}
Let us observe that the equation
%
\begin{equation}
\label{eq:v1} v''(x)-\frac{Q'}{2Q}
v' -Qv = 0 ,
\end{equation}
admits the linearly independent solutions
$v_\pm:= \exp(\pm\xi)$ (the Wronskian of this set of solutions
is $2Q^{1/2}$). We exploit then the variation of constant formula from
which we obtain that
if we find a solution $n(\cdot)$ (say,
in $\bbL^2((-\infty, c])$ for every $c \in\bbR$)
%
\begin{equation}
\label{eq:voc} n(x) = \exp \bigl( \xi(x) \bigr)+ \int_{-\infty}^x
R(y) n(y) \sinh \bigl( \xi (x) -\xi(y) \bigr) \dd y ,
\end{equation}
where $R:= \tilde R/ Q^{1/2}$ (remark that $R$ is bounded and $R(x)=
\mathrm{O}(1/\vert x\vert^4)$ for $\vert x \vert$ large),
then $n(\cdot)$ solves \eqref{eq:n1}. In view of the result, we want
to obtain we set also
%
\begin{equation}
N(x) := n(x) \exp \bigl(-\xi(x) \bigr) ,
\end{equation}
so that \eqref{eq:voc} becomes
%
\begin{equation}
\label{eq:voc2} N(x) = 1+ \int_{-\infty}^x R(y) N(y)
\cK \bigl( \xi(x) -\xi (y) \bigr) \dd y =: 1+\cT_\eta N(x) ,
\end{equation}
where $\cK(z)= (1-\exp(-2z))/2$. We look at $\cT_\eta$ as an
operator that acts on $\bbL^\infty$ functions.

\begin{lemma}
\label{th:Tsmall}
For every $\ga\in(0, \uppi)$ there exists $C>0$ such that, for any
$\eta$ with $\arg(-\eta) \in[-\alpha,\alpha]$,
$\Vert\cT_\eta\Vert_\infty\le C \vert\eta\vert^{-3/4}$.
\end{lemma}

This lemma, applied to \eqref{eq:voc2}, tells us that we can write $N$ as
the operator $\cT_\eta(1-\cT_\eta)^{-1}$\vspace*{1pt} applied to the constant
function equal to $1$
and therefore $\Vert N \Vert_\infty\le2C \vert\eta\vert^{-3/4}$
for $\vert\eta\vert^{3/4} \ge2$.
Therefore, this completes the proof of Theorem~\ref{th:WKBeta}.
\end{pf*}

\begin{pf*}{Proof of Lemma~\ref{th:Tsmall}}
In order to bound $\Vert\cT_\eta\Vert_\infty$, we first need to
look more carefully at the argument of $\cK$ in the integral. By (\ref
{eq:minrealpart}), for any $y<x$ we have
%
\begin{equation}
\label{eq:reK} \Re \bigl(\xi(x) -\xi(y) \bigr) = \int_y^x
\Re \bigl(Q_{\eta}^{1/2}(z) \bigr) \dd z \ge \cos \biggl(
\frac{\alpha+\uppi}{4} \biggr) \sqrt{|\eta|}(x-y) .
\end{equation}
Moreover, by (\ref{eq:maximpart}), we also have for any $y<x$,
%
\begin{equation}
\bigl\vert\Im \bigl(\xi(x) -\xi(y) \bigr) \bigr\vert\le\int_y^x
\bigl\vert \Im \bigl(Q_{\eta}^{1/2}(z) \bigr) \dd z \bigr\vert\le \sqrt{|
\eta|}(x-y).
\end{equation}
Combining the last two inequalities, we further deduce that
%
\begin{eqnarray}
\bigl\vert\Im \bigl(\xi(x) -\xi(y) \bigr) \bigr\vert> \uppi/4 \quad &\Rightarrow& \quad x-y >
\frac{\uppi}{4 \sqrt{|\eta|}} \nonumber\\[-8pt]\\[-8pt]
& \Rightarrow&\quad \Re \bigl(\xi(x) -\xi(y) \bigr) > \frac{\uppi}{4}
\cos \biggl(\frac{\alpha+\uppi}{4} \biggr),\nonumber
\end{eqnarray}
so $ \vert\Im(\xi(x) -\xi(y)) \vert> \uppi/4 $, then $\vert\cK(
\xi(x) -\xi(y))\vert\le(1+ \exp(-\uppi/2)/2)/2< \sqrt{2}/2$.
On the other hand if $\vert\Im(z)\vert\le\uppi/4$, then $\vert\cK
(z) \vert\le\sqrt{2}/2$.
Therefore for any $y<x$,
%
\begin{equation}
\label{eq:boundK} \bigl\llvert \cK \bigl(\xi(x)-\xi(y) \bigr) \bigr\rrvert \le
\frac{\sqrt{2}}2 .
\end{equation}

Furthermore, we see by direct inspection that for every $\alpha\in[0,
\uppi)$ there exists $C_0>0$ such that
%
\begin{equation}
\label{eq:boundR} \bigl\vert R(y) \bigr\vert\le\frac{C_0}{|\eta|+y^4} ,
\end{equation}
for every $y\in\bbR$ and every $\eta$ such that $\arg(-\eta)\in
[-\ga, \ga]$. Therefore
%
\begin{equation}
\label{eq:estT} \Vert\cT_\eta N\Vert_{\infty} \le
C_0 \frac{\sqrt{2}}2 \Vert N\Vert_{\infty} \int
_{-\infty}^\infty \bigl( \vert\eta\vert +y^4
\bigr)^{-1}\dd y = C \eta^{-3/4} \Vert N\Vert_{\infty} ,
\end{equation}
in which the last step is just the definition of $C$.
\end{pf*}

\begin{pf*}{Proof of Theorem~\ref{th:WKBx}} A good deal of this
proof focuses on the asymptotic behavior of $u_{\mathrm{sub},\cG_3}$
as $x \to\infty$.
Recall in fact \eqref{eq:coord}, so that, in view of \eqref{eq:Stransf2},
we have that
%
\begin{equation}
\label{eq:p1WKBx} \lim_{x \to+\infty} u_{\mathrm{sub},\cG_3,\eta}(x) \exp \biggl(-
\frac{x^3}{6} + \mu \frac{x}{2} \biggr) = b(\eta) \lim
_{x \to+\infty} u_{\mathrm{dom},\cG_0, \eta}(x) \exp \biggl(-\frac{x^3}{6} +
\mu \frac{x}{2} \biggr) ,
\end{equation}
provided that the limit on the right-hand side exists.
We actually aim at proving also the analyticity of the left-hand side
in the whole of $\bbC$, but, since
$b(\cdot)$ is entire, this amounts to showing that the limit in the
right-hand side yields an entire function.

Furthermore, recall from Section~\ref{sec:spectrum} that $b(\eta) = 0$ if
and only if $\eta\in\spect$, thus to establish that the left-hand
side of (\ref{eq:p1WKBx}) only vanishes on the spectrum, all we need
to establish is that the limit in the right-hand side does not vanish.

As for what concerns $\eta$ real, we know that $u_{\mathrm{sub},\cG
_3}(x) >0$ for every $x \in\bbR$ if $\eta<\eta_0$ and so the last
statement of Theorem~\ref{th:WKBx} will again follow if we establish
that the limit in the right-hand side is not $0$.

We are going to use an approach parallel to that of the last paragraph
(we are still doing WKB estimates, even if
now $x\to\infty$ instead of $\vert\eta\vert\to\infty$), except
that now, for $\eta\in\R$, $Q_{\eta}$ may take the value zero,
making impossible the change of functions of the last paragraph on the
whole of $\R$.
But now we just working in a neighborhood of $+\infty$ and in fact
the first fact to remark is that when $x$ is large $\Re Q_{\eta}(x)$
is large and $Q_{\eta}(x)$ is in a small sector containing the
positive axis,
that is, $\arg Q_{\eta}(x)$ is small.
We will need also more precise information about $Q_{\eta}(x)$ and we
collect them in the following lemma
for which we introduce the notation $B_r:= \{z\in\bbC\dvt  \vert z \vert
< r\}$:

\begin{lemma}
\label{th:lemt4}
For every $r>0$
there exists $A_0>0$ and $C > 0$ such that for every $A\ge A_0$ we have
(recall that $R$ is defined below \eqref{eq:voc})
%
\begin{equation}
\label{eq:Qest67} \Re \bigl(Q_{\eta}(x) \bigr) \ge\frac{x^4}8,\qquad
\bigl\llvert \arg(Q_{\eta}) (x) \bigr\rrvert \le\frac\uppi4 \quad \mbox{and}\quad
\bigl|R(x)\bigr| \le\frac{C}{x^4} ,
\end{equation}
for $x \ge A$ and
$\eta\in\overline{B}_r$. Moreover
for $x \to\infty$, we have
%
\begin{eqnarray}
\label{eq:taylorQ} %
Q_\eta^{1/2}(x) & =&
\frac{x^2}{2} - \frac{\mu}{2} +\frac{1}{x} +\mathrm{O} \biggl(
\frac{1}{x^2} \biggr) ,\nonumber
\\[-8pt]\\[-8pt]
\frac{1}{Q_\eta^{1/4}(x)} & =& \frac{\sqrt{2}}{x} + \mathrm{O} \biggl(\frac
{1}{x^2}
\biggr) ,\nonumber %
\end{eqnarray}
uniformly for $\eta\in\overline{B}_r$.
\end{lemma}

\begin{pf}
From \eqref{eq:Qeta9} one sees $Q_\eta(x) \sim x^4/4$ for $x \to
\infty$ and that $\Re Q_\eta(x)$ tends to $\infty$,
while $\Im Q_\eta(x)$ stays bounded. Therefore
for any $\gb\in(0, \uppi)$ and any $r>0$ there exists $x_0>0$ such
that $\llvert \arg(Q_{\eta})(x)\rrvert  \le\gb$ for every
$\eta\in\overline{B}_r$ and every $x \ge x_0$.
These observations
yield \eqref{eq:Qest67}: the first two estimates are immediate, the
third requires
making \eqref{eq:voc} explicit.

For \eqref{eq:taylorQ}
we set
$\eta= \eta_1 + \mathrm{i} \eta_2$ and $Q_\eta(x)=r(x)\exp(\mathrm{i} \theta(x))$
so that
%
\begin{equation}
\label{eq:rsqrt0} r(x)^{1/2} = \frac{x^2}{2} \biggl[ \biggl(1-
\frac{2\mu}{x^2}+\frac
{4}{x^3}+\frac{\mu^2 +4\eta_1}{x^4} \biggr)^2 +
\frac{16 \eta
_2^2}{x^8} \biggr]^{1/4} .
\end{equation}
We then observe that for $x$ large (and uniformly in $\eta\in
\overline{B}_r$)
%
\begin{equation}
r(x)^{1/2} = \frac{x^2}{2} -\frac{\mu}{2} +\frac{1}{x}
+ \mathrm{O} \biggl(\frac{1}{x^2} \biggr) \quad \mbox{and}\quad  \theta(x) =
\mathrm{O} \biggl( \frac
{1}{x^{4}} \biggr) ,
\end{equation}
and \eqref{eq:taylorQ} follows.
\end{pf}

In view of the estimates, we are after (or, equivalently, in view of
what the WKB approach suggests)
we now set
%
\begin{equation}
\xi_A (x) := \int_A^x
Q^{1/2}(y) \dd y ,
\end{equation}
and
%
\begin{equation}
N_{\mathrm{dom},A, \eta}(x)=N_{\mathrm{dom},A }(x):= \frac{Q_{\eta}^{1/4}(x) \exp(-\xi_{A}(x)) u_{\mathrm{dom},\cG
_0}(x)}{Q_{\eta}^{1/4}(A)u_{\mathrm{dom},\cG_0}(A)} ,
\end{equation}
for $x \ge A$ (of course we have chosen the normalization so that
$N_{\mathrm{dom},A }(A)=1$).
Of course this requires
%
\begin{equation}
\label{eq:Qeta14} Q_{\eta}^{1/4}(A)u_{\mathrm{dom},\cG_0}(A) \neq0 ,
\end{equation}
but this is granted, uniformly in $\eta\in\overline{B}_r$, for $A$
sufficiently large because $u_{\mathrm{dom},\cG_0}$ is
a dominant solution and because of Lemma~\ref{th:lemt4}. Moreover, by
\eqref{eq:Qest67},
$\eta\mapsto Q_{\eta}^{1/4}(A)$ is entire for $A$ large, while
$u_{\mathrm{dom},\cG_0}(A)$ for every $A$ is entire in $\eta$
as we pointed (see \eqref{eq:subentire} and paragraph leading to
\eqref{eq:coord}), and therefore
the expression in \eqref{eq:Qeta14}
is entire (in $\eta$) for $A$ large.

At this point, we write
%
\begin{equation}
u_{\mathrm{dom},\cG_0, \eta}(x) \exp \biggl(-\frac{x^3}{6} + \mu \frac{x}{2}
\biggr) = \bigl( Q_{\eta}^{1/4}(A)u_{\mathrm{dom},\cG_0}(A) \bigr)
N_{\mathrm{dom},A, \eta}(x) h_\eta(x) ,
\end{equation}
where
%
\begin{equation}
\label{eq:inth} h_\eta(x) := Q_\eta^{-1/4}(x) \exp
\biggl( \int_A^x Q_\eta
^{1/2}(y)\dd y - \frac{x^3}{6} + \frac{\mu}{2}x \biggr) ,
\end{equation}
and we are done if we show that the limits as $x \to\infty$ of
$N_{\mathrm{dom},A, \eta}(x)$ and of $h_\eta(x)$ exist, that they
are non-zero and that
the limit expressions are analytic in $\eta\in B_r$.
Let us notice from now that such a statement for $N_{\mathrm{dom},A,
\eta}(x)$ involves WKB analysis,
while $h_\eta(x)$ is an explicit expression and will be dealt just by
applying the Taylor expansion estimates
\eqref{eq:taylorQ}.

Let us then start with $h_\eta(x)$. For $x$ fixed
$h_{\cdot}(x)$ is now looked upon as a function from $\overline{B}_r$
to $\bbC$ and it is analytic in $B_r$.
By \eqref{eq:taylorQ},
the family of functions
$ \{ h_{\cdot}(x) \}_{x\ge A}$
possesses a limit as $x \to\infty$. Added to that, \eqref
{eq:taylorQ} implies that this family is bounded and bounded away from
$0$ for every $\eta\in\overline{B}_r$, provided
$A$ is sufficiently large. An application of Montel's Theorem \cite{cf:Conway}, page 153, establishes the
analyticity of $\eta\mapsto\lim_{x \to\infty}h_\eta(x)$ in $B_r$.
The fact that the family is bounded away from zero implies of course
that the limit is bounded away from zero.

Let us then turn to $N_{\mathrm{dom},A, \eta}$.
We directly verify exactly like for \eqref{eq:voc2} that $N_{\mathrm
{dom},A, \eta}$ solves for $x \ge A$
%
\begin{equation}
\label{eq:voc2.+} N(x) = 1+ \int_{A}^x R(y) N(y)
\cK \bigl( \xi_{A}(x) -\xi _{A}(y) \bigr) \dd y =: 1+
\cT_{A} N (x) .
\end{equation}
$\cT_A$ is viewed as an operator acting on
$\bbL^\infty([A, \infty);\bbC)$ and one verifies by exploiting
\eqref{eq:Qest67}, exactly like in the proof of Lemma~\ref
{th:Tsmall}, that
$\Vert\cT_A\Vert=\mathrm{o}(1/A^3)$ as $A\to\infty$, uniformly in $\eta\in
\overline{B}_r$ and therefore that \eqref{eq:voc2.+} has a unique
solution $N=N_{\mathrm{dom}, A, \eta}$
satisfying
%
\begin{equation}
\label{eq:boundNdom} \sup_{\eta\in\overline{B}_r}\sup_{x\ge A} \bigl
\llvert N_{\mathrm
{dom}, A, \eta}(x) -1 \bigr\rrvert \le\frac{1}2 ,
\end{equation}
for $A$ sufficiently large. Moreover, by differentiating both sides in
\eqref{eq:voc2.+} (the smoothness of $N_{\mathrm{dom},A}$ follows directly
from the integral equation) we obtain
%
\begin{equation}
\label{eq:Ndomder} N_{\mathrm{dom},A}'(x) = \sqrt{Q_{\eta}(x)}
\int_{A}^x R(y) N_{\mathrm{dom},A}(y) \exp
\bigl(-2 \bigl(\xi_{A}(x)-\xi _{A}(y) \bigr) \bigr) \dd y ,
\end{equation}
so that for every $\eta\in\overline{B}_r$
%
\begin{equation}
\label{eq:Ndomder2} \bigl\vert N_{\mathrm{dom},A}'(x)\bigr\vert\le c
x^2 \exp \bigl(-2 \Re\xi _{A}(x) \bigr) \int
_{A}^x \frac{\exp(2 \Re\xi_{A}(y))}{y^4} \dd y \stackrel{x \to
\infty}\sim\frac{c'}{x^3} ,
\end{equation}
where $c$ and $c'$ are suitable positive constants: we have of course
used the estimates in Lemma~\ref{th:lemt4}, notably \eqref
{eq:taylorQ} and for the asymptotic statement can be obtained by
integration by parts (see, e.g., \cite{cf:BO}, pages 255--256).
Therefore for every $\eta\in\overline{B}_r$, $\lim_{x \to\infty}
N_{\mathrm{dom}, A, \eta}(x) $ exists and
\eqref{eq:Ndomder2} implies also that the limit is non-zero for $A$
sufficiently large. This can be seen also
directly from \eqref{eq:boundNdom}. But \eqref{eq:boundNdom} yields
analyticity too:
since for every $x\ge A$ the function $\eta\mapsto N_{\mathrm{dom},
A, \eta}(x)$ is analytic in $B_r$,
by Montel's Theorem \cite{cf:Conway}, page 153, the limit is analytic in $B_r$.

The proof of Theorem~\ref{th:WKBx} is therefore complete.
\end{pf*}

\subsection{Proof of the corollaries of Section \texorpdfstring{\protect\ref{sec:ODEres}}{4.3}}
\label{sec:proofcorollaries}

\begin{pf*}{Proof of Corollary~\ref{cor:lambda0}}
First of all, recall that the change of variable $\eta=2\gl$, so
$f_\gl= g_{2\gl}$, shows that \eqref{eq:edofgl}
in Corollary~\ref{cor:edo} is the same as \eqref{eq:g}. This is
spelled out also in \eqref{eq:laplace--g}, in a different language.
Therefore, by recalling also \eqref{eq:Stransf} and \eqref{eq:S}, we
see that the properties of the subdominant solution
$u_{\mathrm{sub},\cG_3}$, in particular Remark~\ref{rem:dersub},
yield the existence result assumed in
Corollary~\ref{cor:edo} (and Corollary~\ref{cor:edo2}). Moreover,
Theorem~\ref{th:WKBx} guarantees that the right-hand side of \eqref
{eq:coredo2}, which is actually simply, $1/f_\gl(\infty)$, is
meromorphic, with poles only on the
positive real axis and the first one is $\eta_0/2$. At this point we
use the analytic extension argument detailed in Remark~\ref{rem:analytic}
to get that
\eqref{eq:coredo2} holds for every $\gl< \gl_0=\eta_0/2$.
This establishes \eqref{eq:eta0lambda0}.

For the more precise estimate on this pole claimed in Corollary~\ref
{cor:lambda0}, it suffices to study the behavior
of $g_\eta(\infty)$ near $\eta=\eta_0$, where it vanishes.
Therefore, we set $h_{\eta}(x): = \frac{\mathrm{d} g_{\eta}(x)}{\mathrm{d}\eta}$,
and aim at the existence and
evaluation of $h_{\eta_0}(\infty)$.
We write $h = h_{\eta_0}$ and introduce
%
\begin{equation}
\label{eq:defh} v(x) := \frac{\mathrm{d} u_{\mathrm{sub},\cG_3}}{\mathrm{d}\eta} \biggl\vert_{\eta=
\eta_0} = \exp \biggl(
\frac{x^3}{6} - \mu\frac{x}{2} \biggr) h(x) ,
\end{equation}
so that $v$ is solution to
%
\begin{equation}
\label{eq:v} v'' - Q_{\eta} v = -
u_{\mathrm{sub}, \cG_3} .
\end{equation}
Since we are on the spectrum ($\eta=\eta_0$), that is, $(u_0:=)\,u_{\mathrm
{sub},\cG_3} \propto u_{\mathrm{sub},\cG_0}$, and thus also
%
\begin{equation}
\label{eq:spectrdom} u_{\mathrm{dom}, \cG_3} \propto u_{\mathrm
{dom},\cG_0},
\end{equation}
where we have set $u_{\mathrm{dom}, \cG_3}:= u_{\mathrm{sub }, \cG
_2}$ by using
exactly the same argument that lead to the definition of $u_{\mathrm
{dom},\cG_0}$. Note that we have made the
precise choice (i.e., we choose the multiplicative constant) of
$u_{\mathrm{sub},\cG_3}$ given by \eqref{eq:G3sub}, that we recall here
%
\begin{equation}
\label{eq:u0-} u_{\mathrm{sub},\cG_3} (x) = u_0(x) \stackrel{x\to-
\infty} \sim \exp \biggl(\frac{x^3}{6} -\mu\frac{x}{2} \biggr) .
\end{equation}
By \cite{cf:Sibuya}, Th. 6.1 and Th. 7.2, we have also that there
exists $c>0$ (for the positivity recall that
the solution is positive for $\eta< \eta_0$) such that
%
\begin{equation}
\label{eq:u0+} u_0(x) \stackrel{x\to+\infty} \sim
\frac{c}{x^2}\exp \biggl(-\frac
{x^3}{6} +\mu\frac{x}{2} \biggr) .
\end{equation}
On the other hand, we can also make a precise choice of $U_0 :=
u_{\mathrm{dom},\cG_0}$ by fixing the multiplicative constant in the
asymptotic behavior
of the
subdominant solution in the sector $\cG_1$ (that we use to define the
dominant solution in $\cG_0$, like in Section~\ref{sec:Schr}).
This in turn gives a definite choice of the asymptotic behaviors of $U_0$
%
\begin{equation}
\label{eq:U0} U_0(x) \stackrel{x\to-\infty} \sim
\frac{c'}{x^2} \exp \biggl(-\frac{x^3}{6} +\mu\frac{x}{2} \biggr)
\quad \mbox{and}\quad  U_0(x) \stackrel{x\to\infty} \sim\exp \biggl(
\frac{x^3}{6} -\mu \frac{x}{2} \biggr) ,
\end{equation}
where $c'$ is a constant that appears as a result of fixing to one the
multiplicative constant in the behavior
at $+\infty$ of the dominant solution. Three observations are in order:
\begin{enumerate}[(2)]
\item[(1)] The asymptotic behaviors \eqref{eq:U0} are obtained precisely
like in the proof of Theorem~\ref{eq:WKBx}:
the proof for the case $x \to\infty$ is actually completely contained
in the proof of Theorem~\ref{eq:WKBx}, the case
$x \to-\infty$ requires redoing the asymptotic computation starting
from \eqref{eq:inth}.
\item[(2)] The leading asymptotic behaviors of $u'_0$ and $U'_0$ are
obtained by taking the derivative of the
asymptotic relations \eqref{eq:u0-}, \eqref{eq:u0+} and \eqref
{eq:U0} and by keeping the leading order: this is proven in
\cite{cf:Sibuya}, Th. 6.1, for the subdominant case, while for the
dominant case it just requires a straightforward
generalization of the argument in the proof of Theorem~\ref{th:WKBx},
starting from \eqref{eq:Ndomder} (note that
it is a matter of refining \eqref{eq:Ndomder2}).
\item[(3)] We have made a definite choice of $U_0$ by exploiting the
subdominant solution
in $\cG_1$, which is unique up to a multiplicative constant, and by
fixing the constant
with the second asymptotic statement in \eqref{eq:U0}. Of course we
can replace $U_0$ by adding a multiple of
$u_0$, and this does not change \eqref{eq:U0}: one can clearly see
that also the final result we obtain, that is, \eqref{eq:derfi2},
is invariant under such a change.
\end{enumerate}

Obviously $u_0, U_0$ form a system of independent solutions to (\ref
{eq:S}). It is well known (and straightforward to check) that the
Wronskian $((u_0)'U_0 - (U_0)'u_0)(x)$ is a constant that we call $W$.
$W$ can be computed by using
\eqref{eq:u0-}, \eqref{eq:u0+} and \eqref{eq:U0} along with point
(2) in the list above. Since it can actually be computed
in the two limits $x \to\pm\infty$, one finds $((u_0)'U_0 -
(U_0)'u_0)(\infty)=-c$
and $((u_0)'U_0 - (U_0)'u_0)(-\infty)=c'$, so $c'=-c<0$.
We now use the variation of constants method to give a general expression
\cite{cf:CL}, Th. 6.4,
for $v$ (by Remark~\ref{rem:dersub} $v(-\infty)=0$)
%
\begin{equation}
\label{eq:varconst} v(x) = \frac{1}{W} \biggl(U_0(x) \int
_{-\infty}^{x} u_0^2(y) \dd y
- u_0(x) \int_{-\infty}^x
u_0(y) U_0(y) \biggr) ,
\end{equation}
and for clarity we will substitute $W=-c$ only at the end of the {computation.}
Obviously $u_0 \in\mathbb{L}^2$ so that $\int_{\R} u_0^2(y) \dd y
\in(0,\infty)$.
Also, by \eqref{eq:u0-}, \eqref{eq:u0+} and \eqref{eq:U0}
we see that $u_0 U_0(x) \stackrel{x \to\pm\infty}\sim\frac
{c}{x^2}$, so that also $\int_{\R} u_0 U_0(y) \dd y \in(0, \infty
)$. Therefore,
since $U_0$ is dominant at $+\infty$ whereas $u_0$ is subdominant, it
follows that the second term in the sum of (\ref{eq:varconst}) is
negligible when compared to the first.
Hence,
%
\begin{equation}
v(x) \label{eq:varconst2} \stackrel{x \to\infty} {\sim} \frac{U_0(x)}{W} \int
_{-\infty
}^{\infty} u_0^2(y) \dd y
.
\end{equation}
We now deduce from \eqref{eq:defh} and \eqref{eq:U0} that (here we
insert also $W=-1/c$)
%
\begin{equation}
h(x)\stackrel{x \to\infty} {\longrightarrow} -\frac{1}c \int
_{-\infty}^{\infty} u_0^2(y) \dd y
,
\end{equation}
and therefore
%
\begin{equation}
\label{eq:derfi2} g_{\eta}(\infty) \stackrel{\eta\to\eta_0}
{\sim} (\eta_0-\eta) \int_{-\infty}^{\infty}
\frac{u_0^2(y)} c \dd y .
\end{equation}
Recall that $\Phi_{T_{3,\mu}}(\lambda) = \frac{1}{g_{\eta
}(\infty)}$ and $\eta_0-\eta= 2(\lambda_0-\lambda)$, so we are
done with the proof of Corollary~\ref{cor:lambda0}.
\end{pf*}

\begin{pf*}{Proof of Corollary~\ref{cor:firstorder}}
We recall that the basic formula is \eqref{eq:laplace--g}, that
is Corollary~\ref{cor:edo2}, complemented by Remark~\ref
{rem:analytic} and Remark~\ref{rem:dersub},
and that we found practical at a certain stage to decide that $g_\eta
(-\infty)=1$
(cf. \eqref{eq:Stransf2}). In reality, in this proof the crucial tool
is Theorem~\ref{th:WKBeta},
and that result contains a constant $c$ that we do not determine explicitly.
This actually amounts to saying that for this proof it is easier to
think in terms of
${g_\eta(-\infty)}/{g_\eta(\infty)}$, rather than $1/{g_\eta
(\infty)}$. In any case, by applying Theorem~\ref{th:WKBeta}
one obtains
%
\begin{equation}
\frac{g_\eta(-\infty)}{g_\eta(\infty)} = \lim_{x \to\infty} \frac{g_\eta(-x)}{g_\eta(x)} = \lim
_{x \to\infty} \exp \biggl( -\int_{-x}^{x}
\biggl( Q_\eta^{1/2}(z) -\frac{1}2 z^2 +
\frac\mu2 \biggr) \dd z \biggr) .
\end{equation}
Let us set $\beta_{\eta} = \arg(-\eta)$, so by assumption $\beta
_{\eta} \to\beta$. We are aiming at proving that if we fix any $\gb
_0\in(0, \uppi)$, for every $\gd>0$ there exists $\kappa>0$ such that
%
\begin{eqnarray}
\label{eq:firstorder} &&\biggl\llvert \frac{1}{|\eta|^{3/4}} \lim_{x \to\infty}
\Re \biggl[ \int_{-x}^{x} \biggl(
Q_\eta^{1/2}(z) -\frac{1}2 z^2 + \frac
\mu2 \biggr) \dd z \biggr] - c (\beta_\eta ) \biggr\rrvert \le\gd\nonumber\\[-8pt]\\[-8pt]
&&\quad  \mbox{for } \vert\eta\vert> \kappa\mbox{ and } \vert\gb_\eta\vert\le
\gb_0 ,\nonumber
\end{eqnarray}
where $c(\beta_{\eta})$ is as in \eqref{eq:a_0}.

In order to establish this,
recall that $Q_{\eta}(z) = z+\frac{1}{4}(z^2-\mu)^2 - \eta$, so
that the leading order depends on wether $|z|^4$ or $|\eta|$ is the
largest. More precisely, there exists a $A_0$ large enough so that, for
any $z$ such that $|z| \ge A |\eta|^{1/4}$ with $A \ge A_0$, we find
(similarly to (\ref{eq:taylorQ}))
%
\begin{equation}
\label{eq:taylorQeta} \biggl\llvert Q_{\eta}^{1/2}(z) -
\frac{z^2}{2} + \frac{\mu}{2} -\frac{1}{z} \biggr\rrvert \le C
\frac{|\eta|}{z^2},
\end{equation}
where $C$ is a positive constant which does not depend on $A$ nor $\eta$.
Using the fact that $x \mapsto1/x$ is odd, and then (\ref
{eq:taylorQeta}), we find that
%
\begin{eqnarray}
\label{eq:neglaboveA}
&&\lim_{x \to\infty} \biggl\llvert \int
_{x \ge|z| \ge A |\eta|^{1/4}} \biggl( Q_{\eta}^{1/2}(z) -
\frac{z^2}{2} + \frac{\mu}{2} \biggr) \dd z \biggr\rrvert\nonumber
\\
&&\qquad = \lim_{x \to\infty} \biggl\llvert \int_{x \ge|z| \ge A |\eta|^{1/4}}
\biggl( Q_{\eta}^{1/2}(z) - \frac{z^2}{2} +
\frac{\mu}{2} -\frac
{1}{z} \biggr) \dd z \biggr\rrvert
\\
&&\qquad \le\lim_{x \to\infty} \int_{x
\ge|z|\ge A |\eta|^{1/4}}
\frac{C |\eta|}{z^2} = 2\frac{C}{A}|\eta|^{3/4}.\nonumber
\end{eqnarray}
On the other hand,
%
\begin{eqnarray}
&& \int_{-A |\eta|^{1/4}}^{A |\eta|^{1/4}} \biggl( Q_{\eta}^{1/2}(z)
-\frac{1}2 z^2 + \frac\mu2 \biggr) \dd z
\nonumber\\[-8pt]\\[-8pt]\nonumber
&&\qquad = |\eta|^{1/4} \int_{-A}^{A} \biggl(
Q_{\eta}^{1/2} \bigl(y|\eta|^{1/4} \bigr) -
\frac{|\eta|^{1/2} y^2}{2} + \frac\mu2 \biggr) \dd y .
\end{eqnarray}
It remains to estimate the leading order, uniformly for $|\eta| $
large and in the sector $\vert\gb_\eta\vert\le\gb_0$
(of course we can assume $\gb_0>\uppi/2$), of $\Re(Q_{\eta
}^{1/2}(y|\eta|^{1/4}))$.
Once again we set
$Q_{\eta}(x)=r_{\eta}(x)\exp(\mathrm{i} \theta(x))$ and
\begin{itemize}
\item We find that the norm of $Q_{\eta}^{1/2}(y|\eta|^{1/4})$ is
%
\begin{equation}
\label{eq:uclaim1} r_{\eta}^{1/2} \bigl(|\eta|^{1/4} y
\bigr) = |\eta|^{1/2} \biggl[1+\cos (\beta ) \frac{y^4}{2}+
\frac{y^8}{16} \biggr]^{1/4} + \mathrm{o} \bigl(|\eta
|^{1/2} \bigr) .
\end{equation}
\item The cosine of the argument of $Q_{\eta}^{1/2}(y|\eta|^{1/4})$
is (again we use $\beta_{\eta} \to\beta$)
%
\begin{eqnarray}
\label{eq:uclaim2} \cos \biggl(\frac{1}{2}\theta \bigl(|
\eta|^{1/4} y \bigr) \biggr) & =& \frac{1}{\sqrt{2}} \sqrt{1+
{1}\Bigl/{\sqrt{1+\frac{\sin
^2(\beta)}{\cos^2(\beta)+\cos(\beta)\sfrac{y^4}{2} + \sfrac
{y^8}{16}}} }} + \mathrm{o}(1)\nonumber
\\[-8pt]\\[-8pt]
& =& \frac{1}{\sqrt{2}} \sqrt{1+ \frac{\cos(\beta) +\sfrac
{y^4}{4}}{\sqrt{1 +\cos(\beta)\sfrac{y^4}{2} +\sfrac{y^8}{16}}}} + \mathrm{o}(1) ,\nonumber
\end{eqnarray}
where $\mathrm{o}(1)$ is as $\vert\eta\vert\to\infty$ in the sector $\vert
\gb\vert\le\uppi/2$ (i.e., $\cos(\beta) \ge0$, see below for the other
cases).
\end{itemize}

The estimate \eqref{eq:uclaim1} above is a simple consequence of the
exact expression for $r_{\eta}^{1/2}(x)$, cf. \eqref{eq:rsqrt0}, and
the fact that $\beta_{\eta} \to\beta$.

As for \eqref{eq:uclaim2}, we of course also use that $\beta_{\eta}
\to\beta$, but then recall we made the extra assumption $\cos(\beta
) \ge0$ and write
%
\begin{equation}
\theta \bigl(|\eta|^{1/4} y \bigr) = \arctan \biggl(\frac{\sin
(\beta)}{\cos(\beta)+\sfrac{y^4}{4}}
\biggr) + \mathrm{o}(1) ,
\end{equation}
with $\mathrm{o}(1)$ as in \eqref{eq:uclaim2},
and then simply use the fact that
%
\begin{equation}
\cos \biggl(\frac{1}{2}\arctan(x) \biggr) = \frac{1}{\sqrt{2}} \sqrt{1+
\frac{1}{\sqrt{1+x^2}}}.
\end{equation}
Combining \eqref{eq:uclaim1} and \eqref{eq:uclaim2} yields, in the
case $\cos(\beta) \ge0$,
%
\begin{eqnarray}
&&\hspace*{-20pt}\Re \bigl(Q_{\eta}^{1/2} \bigl(y|
\eta|^{1/4} \bigr) \bigr) \nonumber\\
&&\hspace*{-20pt}\quad = |\eta|^{1/2} \biggl[1+\cos(\beta)
\frac{y^4}{2}+ \frac
{y^8}{16} \biggr]^{1/4}\sqrt{1+
\frac{\cos(\beta) +\sfrac
{y^4}{4}}{\sqrt{1 +\cos(\beta)\sfrac{y^4}{2} +\sfrac{y^8}{16}}}} + \mathrm{o} \bigl(| \eta|^{1/2} \bigr)
\\
&&\hspace*{-20pt}\quad = |\eta|^{1/2} \frac{1}{\sqrt{2}} \sqrt{\sqrt{1 +\cos(\beta)
\frac
{y^4}{2} +\frac{y^8}{16}} + \cos(\beta) +\frac{y^4}{4}} +
\mathrm{o} \bigl(|\eta |^{1/2} \bigr) ,\nonumber %
\end{eqnarray}
as long as $\cos(\gb)\ge0$.
The cases $\cos(\beta) < 0, \sin(\beta) >0$ and $\cos(\beta)<0,
\sin(\beta)<0$ are treated in a similar way, but then in the
expression of the argument we need to add or subtract $\uppi$ to
$\arctan (\frac{\sin(\beta)}{\cos(\beta)+\sfrac
{y^4}{4}} )$, and then use a similar expression for $\sin(\frac
{1}{2}\arctan(x))$. In addition for these estimates one needs to use
the fact that $\vert\beta_{\eta} \vert\le\gb_0
\in(\uppi/2, \uppi)$, so that $\cos(\gb_{\eta})$ is bounded away from $-1$.
Elementary algebraic manipulation then lead to the exact same
expression for $\Re(Q_{\eta}^{1/2}(y|\eta|^{1/4}))$ as above.

We have therefore found
%
\begin{eqnarray}
&&\Re \biggl[ |\eta|^{1/4} \int_{-A}^{A}
\biggl( Q_{\eta}^{1/2} \bigl(y|\eta |^{1/4} \bigr) -
\frac{|\eta|^{1/2} y^2}{2} + \frac{\mu}{2} \biggr) \dd y \biggr]
\nonumber\\[-8pt]\\[-8pt]
&&\quad = |\eta|^{3/4} \int_{-A}^A
f_0^{(\beta)}(y) \dd y + \mathrm{o} \bigl(|
\eta|^{3/4} \bigr) ,\nonumber
\end{eqnarray}
uniformly when $\vert\gb_{\eta} \vert\le\gb_0$.
Along with (\ref{eq:neglaboveA}), this proves \eqref{eq:firstorder}
and therefore the proof of Corollary~\ref{cor:firstorder} is complete.
\end{pf*}

\begin{pf*}{Proof of Corollary~\ref{cor:sharp}}
Let us now turn to the more detailed estimate of Corollary~\ref
{cor:sharp} in the case of $\gr:=-\eta>0$.
By the same change of variables $z=y \eta^{1/4}$ we get
%
\begin{equation}
\frac{g_\eta(-\infty)}{g_\eta(\infty)} = \exp \biggl( -\gr ^{1/4} \int
_{-\infty}^{\infty} \biggl( Q_{\eta}^{1/2}
\bigl(y \gr^{1/4} \bigr) - \frac{\gr^{1/2} y^2}{2} + \frac\mu2 \biggr) \dd y
\biggr) .
\end{equation}
In this case elementary computations lead to
%
\begin{equation}
\gr^{1/4} \biggl(\sqrt{Q_{\eta} \bigl(\gr^{1/4}y
\bigr)} -\frac{1}2 \gr^{1/2}y^2 + \frac\mu2 \biggr) =
\gr^{3/4} f_0(y) + \gr^{1/4} f_1(y) +
\frac{y}{\sqrt{y^4+4}}+ r_\eta(y) ,
\end{equation}
where $\vert r_\eta(y) \vert\le C/(\gr^{1/4}(1+ y^2))$.
\end{pf*}

\section{ODE analysis: The pitchfork case}
\label{sec:ode2}

We consider here the even $d$ case. Like in the previous section
we set $\eta=2\lambda\in\mathbb{C}$, but now we look at even
solutions to the equation
$g''-V_{\mu}'g'+\eta g=0$, that is, we look for the (unique) solution
$g=g_\eta$ to
%
\begin{equation}
\label{eq:geven} g''(x) + V_{\mu}'(x)
g'(x) + \eta g(x) = 0,\qquad  g(0)=1 , g'(0)=0 .
\end{equation}
Again we set $g(x) = u(x) \exp(V_{\mu}(x)), x \in\R$, which leads
to the unique solution to
%
\begin{equation}
\label{eq:Squart} u''(x) - Q_{\eta}(x)u(x) =
0 , \qquad u(0)=1, u'(0)=0 ,
\end{equation}
with, like in \eqref{eq:Qeta9},
$Q_{\eta}(x) = (V_{\mu}')^2(x)-V_{\mu}''(x) - \eta$.
Since $V_\mu(\cdot)$ is even, the function $Q_{\eta}(\cdot
)$ is even too. Therefore when $v(\cdot)$ is \emph{any} given
solution to $u'' - Q_{\eta}u = 0$, also $v(-\cdot)$ is a solution.
In addition when $v$ \emph{does not vanish at the origin}, then $u =
(v(\cdot)+v(-\cdot))/2v(0)$ is also a solution to $u'' - Q_{\eta}u =
0$ and in fact $u$ is the unique
solution to (\ref{eq:Squart}) as it obviously satisfies $u(0)=1, u'(0)=0$.
Since our solution to (\ref{eq:geven}) is $g=u \exp(V_{\mu})$ it is
also even.

This simple argument can be applied to a solution to $u'' - Q_{\eta}u
= 0$ which is subdominant in a given sector.
More precisely Theorem~6.1 in \cite{cf:Sibuya} guarantees existence of
subdominant (and dominant) solutions, and they are defined in the open
sectors $\cG_k = \{x \in\bbC\dvt  \llvert \arg(x) - \frac{k\uppi
}{d}\rrvert  < \frac{\uppi}{2d}\}$, $k \in\{-d+1,-d+2,\ldots,d-1,d\}$
(see Figure~\ref{fig:sectors}). Of course $\cG_d$ contains $(-\infty
,0)$ while $\cG_0$ contains $(0, \infty)$. More interestingly,
because of the symmetry with respect to the origin we know that
$u_{\mathrm{sub}, \cG_d}(\cdot) = u_{\mathrm{sub}, \cG_0}(-\cdot)$.

Letting
%
\begin{equation}
u_1 := c_1 (u_{\mathrm{sub},\cG_0}+u_{\mathrm
{sub},\cG_d}), \qquad \mbox{with } c_1 := %
\cases{ 1/ \bigl(2 u_{\mathrm{sub},\cG_0}(0)
\bigr) & \quad \mbox{if} $u_{\mathrm{sub},\cG_0}(0) \ne0$,
\cr
1 & \quad \mbox{otherwise}, }
\end{equation}
we can use our previous reasoning to see that when $u_{\mathrm
{sub},\cG_0}(0) \ne0$ (as will see this is always the case for $\eta
$ left of the spectrum),
then $u_1$ must be the solution to (\ref{eq:Squart}).

In the second case $u_{\mathrm{sub},\cG_0}(0)=0$, then $u_1$ is the
solution to $u'' - Q_{\eta}u = 0, u(0)=u'(0)=0$ so that in fact $u_1
\equiv0$.

Similarly, we let
%
\begin{equation}
u_2 := c_2 (u_{\mathrm{dom},\cG_0}+u_{\mathrm
{dom},\cG_d}), \qquad \mbox{with } c_2 := %
\cases{ 1/ \bigl(2 u_{\mathrm{dom},\cG_0}(0)
\bigr) & \quad \mbox{if} $u_{\mathrm{dom},\cG_0}(0) \ne0$,
\cr
1 & \quad \mbox{otherwise}, }
\end{equation}
where, in strict analogy with what we have done in Section~\ref{sec:Schr},
we made the choice
$u_{\mathrm{dom},\cG_0}:= u_{\mathrm{sub},\cG_1}$, and by symmetry again
$u_{\mathrm{dom},\cG_0} (\cdot) = u_{\mathrm{dom},\cG_d}(-\cdot)$.
Again by the same reasoning, when $u_{\mathrm{dom},\cG_0}(0) \ne0$
(as we will see this always happens when $\eta$ is strictly left of
the spectrum), then $u_2$ is the solution to (\ref{eq:Squart}),
whereas in the case $u_{\mathrm{dom},\cG_0}(0) = 0$, we find $u_2
\equiv0$.

Note finally that $u_{\mathrm{sub},\cG_0}$ and $u_{\mathrm{dom},\cG_0}$
form a basis, thus $u_1\equiv0$ and $u_2\equiv0$ cannot happen at the
same time and $u_\eta=u=(u_1+u_2)/(u_1(0)+u_2(0))$ always is the
solution to
\eqref{eq:Squart}.

In \cite{cf:CL}, Ch. 9, Problem 1, one can find a proof of the fact
that \eqref{eq:Squart}
admits a $\bbL^2([0,\infty))$ solution if and only if $\eta\in
\spect=\{\eta_0, \eta_1, \ldots\}$, with $\eta_j< \eta_{j+1}$.

Such solution has to be proportional to $u_{\mathrm{sub},\cG_0}$
because it has a zero limit at infinity. In particular $u_{\mathrm
{sub}, \cG_0}(0) \ne0$. By unicity of the solution to (\ref
{eq:Squart}), we conclude (as in the odd case) that when $\eta\in
\spect$, $u_{\mathrm{sub}, \cG_0} = u_{\mathrm{sub},\cG_d}$ and we
have in fact a $\bbL^2((-\infty,\infty))$ solution.

In \cite{cf:CL}, Ch. 9, Problem 1, there is also a precise
characterization of the number and locations of the zeros of the
solutions to \eqref{eq:Squart} and, in particular,
for $\eta\le\eta_0$ the solution to \eqref{eq:Squart} does not
change sign, that is, $u(\cdot)>0$. As announced earlier, $u_{\mathrm
{sub}, \cG_0}(0) \ne0$ as long as $\eta\le\eta_0$.

Because of Theorem~7.1 in \cite{cf:Sibuya}, also $u_{\mathrm{dom},
\cG_0}(0) \ne0$ as long as $\eta< \eta_0$.
Moreover, notice that it is exactly when $\eta\in\spect$ that we
have $u_2\equiv0$, that is $u_{\mathrm{dom},\cG_0}$ is odd.

Regularity properties turn out to be easier than in the odd case.
Indeed note that we have chosen to fix $g(0)$, and therefore $u(0)$, to
$1$, along with the zero slope condition. In the
odd $d$ case instead we had to put a boundary condition at $-\infty$,
$g(-\infty)=1$, which contains
information both on the slope and the size of the function: dealing
with these \textsl{non-standard} boundary conditions has required the
approach developed in \cite{cf:Sibuya}, notably for the regularity
issues. It is clear that the boundary conditions in the
even $d$ case are more standard and
results like the analytic dependence of solutions on the parameter of
the equations are also standard, see, for example, \cite{cf:CL}, Ch. 3.

Nevertheless the fact that $u_\eta(x)$ is entire both in $\eta$ and
$x$ can
also be directly extracted from the representation formula we have just
obtained in terms of subdominant and dominant solutions,
of course by exploiting the analyticity properties of (sub)dominant
solutions (this provides an approach alternative to \cite{cf:CL}, Ch. 3).


The considerations we just made directly imply that for the solution
$u=u_\eta$ of \eqref{eq:Squart} one can write
%
\begin{equation}
\label{eq:standard-d} u = a u_{\mathrm{sub},\cG_0} + b u_{\mathrm{dom},\cG_0} ,
\end{equation}
where $a=a(\eta)$, $b=b(\eta)$ and $a(\cdot)$ and $b(\cdot)$ are
entire functions. Moreover, as in the odd case, $b(\eta)=0$ if
and only if $\eta\in\spect$.

Of course the sharp behavior for $x$ large of $u_{\mathrm{sub},\cG
_0}$ can be taken
from \cite{cf:Sibuya}, Th. 6.1, and the one of $u_{\mathrm{dom},\cG
_0}$ can be derived by
the argument in the proof of Theorem~\ref{th:WKBeta}.
As usual, these solutions are defined
up to a multiplicative constant.

By the exact same techniques as in Section~\ref{sec:ode}, and similar
to Theorems \ref{th:WKBeta}, \ref{th:WKBx}, we are further able to
describe precisely (up to a multiplicative constant) the asymptotic
behaviour of $u_{\eta}$, both in the limit $\eta\to-\infty$, and $x
\to\infty$. As before we define $Q_{\eta}^{1/2}(x)$ (resp., $Q_{\eta
}^{1/4}$(x)) the square root of $Q_{\eta}(x)$ (resp., of $Q_{\eta
}^{1/2}(x)$) satisfying $\arg(Q_{\eta}^{1/2}(x)) \in(-\uppi/2, \uppi
/2]$ (resp., $\arg(Q_{\eta}^{1/4}(x)) \in(-\uppi/4, \uppi/4]$).

\begin{theorem}
\label{th:WKBetaquartic}
Fix $\mu\in\bbR$ and $\alpha\in(0,\uppi)$. For any $\eta\in
\mathbb{C} \setminus(0,\infty)$ there exists a solution $u$ to
$u''-Q_{\eta}u=0$ satisfying, uniformly in $\eta\in\mathbb{C}$ such
that $\vert\arg(\eta) \vert< \alpha$,
%
\begin{equation}
\label{eq:WKBetaquartic} \sup_{x \in\bbR} \biggl\llvert u(x)
Q_\eta^{1/4}(x) \exp \biggl(-\int_0^x
Q_\eta^{1/2} (y) \dd y \biggr) -1 \biggr\rrvert = \mathrm{O}
\bigl( | \eta|^{-d/(d+1)} \bigr) .
\end{equation}
\end{theorem}

\begin{theorem} \label{th:WKBxquartic}
For $\eta\in\mathbb{C} $ the following limit exists:
%
\begin{equation}
\label{eq:WKBxquartic} \lim_{x \to\infty} u_{\mathrm{sub},\cG_d}(x) 
\exp \bigl(V_\mu(x) \bigr) =: \ell(\eta) .
\end{equation}
Moreover $\ell(\cdot)$ is entire, $\ell(\eta)=0$ if and only if
$\eta\in\spect$ and $\ell(\eta)>0$ for $\eta< \eta_0=\min
(\spect)$.
\end{theorem}

We now turn to the asymptotic analysis of the Laplace transform when
$\eta\nearrow\eta_0$, when $|\eta| \to\infty$ along a ray which
is not the positive half line, and the more precise result in the
particular case when $\eta\to-\infty$ (the latter will only be
stated in the case $d=4$).

Recall that $\gl_0(T_{d, \mu})$ is defined in \eqref{eq:lambda0} and
\eqref{eq:lim}.

\begin{cor} \label{cor:lambda0quart}
We have that $2\gl_0(T_{d, \mu})=\eta_0$. Moreover, $\Phi_{T_{d,
\mu}}(\cdot)$ extends as a meromorphic function to $\bbC$ and
it has a simple pole in $\gl_0$ with residue
%
\begin{equation}
\label{eq:residue4} C_{d, \mu} = -\frac{u_{\mathrm{sub},\cG_0, \eta_0}(0)}{\int_\bbR
u^2_{\mathrm{sub},\cG_0, \eta_0}(x)\dd x } .
\end{equation}
\end{cor}

The proof is very similar to the proof of Corollary~\ref{cor:lambda0}.
It goes through the differentiation step
\eqref{eq:v} at $\eta=\eta_0$.

In the even case,
we have $u_{0} := u_{\mathrm{sub},\cG_0, \eta_0}=u_{\mathrm
{sub},\cG_d, \eta_0}$ by symmetry, where we choose the multiplicative
factor so that
%
\begin{equation}
\label{eq:sub6} u_{0}(x) \stackrel{x \to\pm\infty}\sim
\frac{1}{\vert x \vert
^{d-1}} \exp \bigl( V_{\mu}(x) \bigr) .
\end{equation}
Equation \eqref{eq:sub6} simplifies the computation of the Wronskian
that this time can be chosen equal to one (by properly choosing the
dominant solution $U_0$ so that $U_0 \sim_{x \to-\infty} \frac
{1}{2}\exp(-V_{\mu}(x))$), this yields
$ h(x) = v(x)\exp(V_{\mu}(x)) \sim_{x \to\infty} \frac{-1}{2}
\int_{\bbR} u_{0}^2(x) \,\mathrm{d}x$.

It follows that $g_{\eta_0}(\infty)\sim(\eta_0-\eta)
\int_\bbR u_0^2(x) \dd x$, but
we shall not forget that in the even case $\Phi_{T_{d,\mu}}(\lambda)
= \frac{g_{\eta}(0)}{g_{\eta}(\infty)}$, and from this (and the
fact that $\eta_0-\eta= 2(\lambda_0-\lambda)$), \eqref
{eq:residue4} follows.

When $|\eta| \to\infty$ on a ray that is not the positive half-line,
we are able to deduce the following result, analogous to Corollary~\ref
{cor:firstorder}. Recall here that $\beta= \arg(-\eta)$.

\begin{cor}
\label{cor:firstorder-quart}
Fix $\beta_0 \in(0,\uppi)$, uniformly in $\gb\in[-\gb_0, \gb_0]$
we have that
%
\begin{equation}
\label{eq:firstorder0-quart} \limtwo{|\eta| \to\infty} { \arg(-\eta) \to\beta} \vert\eta
\vert^{\afrac{-d}{2(d-1)}}\log g_\eta(\infty) = -c(\gb) ,
\end{equation}
%
where $c(\beta) := \int_{0}^{\infty} f_0^{(\beta)}(y) \dd y$, and
%
\begin{equation}
\label{eq:a_0-quart} f_0^{(\beta)}(y) = \frac{1}{\sqrt{2}}
\sqrt{y^{2(d-1)}+\cos(\beta ) + \sqrt{1+2\cos(\beta) y^{2(d-1)} +
y^{4(d-1)}}} - y^{d-1}.
\end{equation}
\end{cor}


We now turn to the particular case $d=4$ to write the details of the
precise asymptotics when $\eta\to-\infty$ (general even $d$ can be
treated with a similar method, but the general expressions are cumbersome).
Note that in this case we are looking at
%
\begin{equation}
V_{\mu}(x) = -\frac{x^{4}}{4} + \mu\frac{x^2}{2},\qquad
Q_{\eta}(x) = \bigl(x^3-\mu x \bigr)^2 + 3
x^2 - \mu- \eta.
\end{equation}

\begin{cor} \label{cor:etareal-quart}
Consider the solution $g=g_\eta$ to \eqref{eq:geven}.
As $\eta\to-\infty$ we have
%
\begin{equation}
\label{eq:etareal-quart} \frac{1}{g_\eta(\infty)} = 2|\eta|^{1/4} \exp \bigl( - \vert
\eta \vert^{2/3} F_0 - \vert\eta\vert^{1/3}
F_1 - F_2 + \mathrm{O} \bigl({\vert \eta
\vert^{-1/3}} \bigr) \bigr),
\end{equation}
where for $i=0,1,2$, $F_i = \int_{0}^{\infty} f_i(y) \,\mathrm{d}y$, and
%
\begin{eqnarray}
f_0(y) &=& \sqrt{1+y^6} - y^3
,\nonumber
\\
f_1(y) &=& \mu \biggl(y-\frac{y^4}{\sqrt{y^6+1}} \biggr) ,
\nonumber\\[-8pt]\\[-8pt]
f_2(y) &=& \frac{1}{2} \biggl( \frac{(\mu^2+3)y^2}{\sqrt{y^6+1}} -
\frac{3}{y+1}-\frac{\mu^2 y^8}{(y^6+1)^{3/2}} \biggr)
\nonumber\\
&=& \frac{1}2\frac{\mathrm{d}}{\mathrm{d} y} \biggl( \bigl( \mathrm{arcsinh}
\bigl(y^3 \bigr) - 3\mu^2\log(1 + y) \bigr) +
\frac{y^3}{
3 \sqrt{1 + y^6}
} \biggr) .\nonumber %
\end{eqnarray}
\end{cor}

Note that $F_2= \frac{1}2\log2 + \frac{\mu^2}6$. $F_1$ and $F_2$ can
be written in terms of \textsl{elliptic integrals of the first
and second kind}:
for $\phi\in(-\uppi/2, \uppi/2)$ and $m <1$
%
\begin{equation}
\label{eq:Ephim} \mathtt{E}_\pm(\phi\vert m) := \int
_0^\phi \bigl(\sqrt{1 - m \sin^2 (
\theta)} \bigr)^{\pm1} \dd\theta.
\end{equation}
With this definition and by setting $q_\pm= \sqrt{3} \pm2$, we have
%
\begin{equation}
\label{eq:kappa0} F_0 = \frac{3^{3/4}}8 \mathtt{E}_- \bigl(
\mathrm{arccos} (q_- ) \vert q_+/4 \bigr) ,
\end{equation}
and
%
\begin{eqnarray}
\label{eq:kappa1} F_1 &=& \frac\mu{12}\bigl( -3 + 3 \sqrt{3} \nonumber
\\[-8pt]\\[-8pt]
&&\hphantom{\frac\mu{12}\bigl(}{} + 3^{1/4} \bigl(6 \mathtt{E}_+ \bigl(
\mathrm{arccos}(q_-) \vert q_+/4 \bigr) + (-3 + \sqrt{3}) \mathtt{E}_- \bigl(
\mathrm{arccos}(q_-) \vert q_+/4 \bigr) \bigr) \bigr) .\nonumber
\end{eqnarray}
From this, we derive the constants that appear in the statement
of Theorem~\ref{th:Phileft}:
%
\begin{equation}
\label{eq:forPl} C_{2/3}:=2^{2/3}F_0 = 1.6693\!
\ldots\quad \mbox{and}\quad  C_{1/3} =\frac
{2^{1/3 } F_1} \mu= 0.5432\!\ldots.
\end{equation}

%


\section*{Acknowledgements}
We thank K. Pakdaman for having pointed out to us \cite{cf:SH}.
G.G. acknowledges the support of ANR Grants SHEPI and ManDy.



\printhistory
\end{document}